\documentclass[12pt,oneside]{book}
\usepackage{geometry}                		
\usepackage{graphicx}									
\usepackage{amssymb}\vspace{0.25cm}
\usepackage{mathtools}
\usepackage{setspace}
\usepackage{tikz-cd}
\usepackage[mathscr]{euscript}
\newcommand{\abs}[1]{\left\vert#1\right\vert}
\usepackage[normalem]{ulem}
\usepackage{manfnt}
\usepackage{pgfplots}
\usepackage{pifont}
\usepackage[safe]{tipa}
\usepackage{marvosym}
\usepackage{scalerel,stackengine}
\usepackage{titlesec}
\usepackage{makeidx}
\usepackage{centernot}
\usepackage{xspace}
\usepackage[noauto]{chappg}
\usepackage[OT2,T1]{fontenc}
\usepackage{hyperref}
\usepackage{amsmath}
\usepackage[amsmath,thmmarks]{ntheorem}
\usepackage{tabularx}
\usepackage{yfonts}
\usepackage{epstopdf}
\usepackage{pdfpages}
\newcolumntype{Y}{>{\raggedright\arraybackslash}X}

\newcommand\N{%
\mathbb{N}
}

\newcommand\R{%
\mathbb{R}
}

\newcommand\X{%
\mathbb{X}
}

\newcommand\Z{%
\mathbb{Z}
}

\newcommand\sB{%
\mathcal{B}
}

\newcommand\sD{%
\mathcal{D}
}
\newcommand\sE{%
\mathcal{E}
}
\newcommand\sF{%
\mathcal{F}
}
\newcommand\sH{%
\mathcal{H}
}
\newcommand\sI{%
\mathcal{I}
}

\newcommand\sP{%
\mathcal{P}
}
\newcommand\sR{%
\mathcal{R}
}

\newcommand\tBC{%
\text{BC}
}
\newcommand\tqa{%
\text{qa}
}
\newcommand\tqsa{%
\text{qsa-}
}
\newcommand\tBO{%
\text{BO}
}

\newcommand\td{%
\text{d}
}

\newcommand\tB{%
\text{B}
}

\newcommand\Lp{
\text{L}
}

\newcommand\tU{%
\text{U}
}

\newcommand\tT{%
\text{T}
}

\def\thereforex{\boldsymbol{\text{ }
\leavevmode
\lower0.4ex\hbox{\textbullet}
\kern-.9em\raise1.1ex\hbox{\textbullet}
\kern-0.9em\lower0.4ex\hbox{\textbullet}
\hspace{0.1cm}\thinspace\text{ }}}
\def\thereforez{\boldsymbol{\text{ }
\leavevmode
\lower0.4ex\hbox{$\circ$}
\kern-.9em\raise1.1ex\hbox{$\circ$}
\kern-0.9em\lower0.4ex\hbox{$\circ$}
\hspace{0.1cm}\thinspace\text{ }}}
\newcommand\ra{%
\rightarrow
}

\newcommand\ds{%
\displaystyle
}

\newcommand\diam{%
\text{diam}\hspace{0.05cm}
}

\newcommand\un[1]{%
\underline{#1}\xspace
}

\newcommand\hsx{%
\hspace{0.05cm}
}
\newcommand\hsy{%
\hspace{0.03cm}
}

\newcommand{\norm}[1]{\left\lVert #1 \right\rVert}

\newcommand{\normx}[1]{\big|\hspace{-.05cm}\big| #1 \big|\hspace{-.05cm}\big|}
\newcommand{\normxx}[1]{\bigg|\hspace{-.05cm}\bigg| #1 \bigg|\hspace{-.05cm}\bigg|}

\newcommand\BV{%
 \text{BV}
}

\newcommand\AC{%
\text{AC}
}

\definecolor{ultramarine}{RGB}{0, 32, 96}

\definecolor{darkcerulean}{rgb}{0.3, 0.27, 0.49}

\definecolor{forestgreen}{rgb}{0.0, 0.27, 0.13}
\definecolor{forestgreenweb}{rgb}{0.13, 0.55, 0.13}
\definecolor{deepjunglegreen}{rgb}{0.0, 0.29, 0.29}

\definecolor{midnightblue}{rgb}{0.1, 0.1, 0.44}
\definecolor{midnightgreen}{rgb}{0.0, 0.29, 0.33}
\definecolor{myrtle}{rgb}{0.13, 0.26, 0.12}
\definecolor{darkviolet}{rgb}{0.58, 0.0, 0.83}
\definecolor{darkgreen}{rgb}{0.0, 0.2, 0.13}
\definecolor{officegreen}{rgb}{0.0, 0.5, 0.0}

\definecolor{harvardcrimson}{rgb}{0.79, 0.0, 0.09}
\definecolor{hollywoodcerise}{rgb}{0.96, 0.0, 0.63}
\definecolor{debianred}{rgb}{0.84, 0.04, 0.33}
\definecolor{darkturquoise}{rgb}{0.0, 0.81, 0.82}

\definecolor{darktangernine}{rgb}{1.0, 0.66, 0.07}
\definecolor{aureolin}{rgb}{0.99, 0.93, 0.0}
\definecolor{canaryyellow}{rgb}{1.0, 0.94, 0.0}
\definecolor{amber}{rgb}{1.0, 0.75, 0.0}

\definecolor{urobilin}{rgb}{0.88, 0.68, 0.13}
\definecolor{uscgold}{rgb}{1.0, 0.8, 0.0}

\usepackage{scalerel,stackengine}
\stackMath
\newcommand\reallywidehat[1]{%
\savestack{\tmpbox}{\stretchto{%
  \scaleto{%
    \scalerel*[\widthof{\ensuremath{#1}}]{\kern-.6pt\bigwedge\kern-.6pt}%
    {\rule[-\textheight/2]{1ex}{\textheight}}
  }{\textheight}%
}{0.5ex}}%
\stackon[1pt]{#1}{\tmpbox}%
}
\makeatletter
\DeclareRobustCommand\widecheck[1]{{\mathpalette\@widecheck{#1}}}
\def\@widecheck#1#2{%
    \setbox\z@\hbox{\m@th$#1#2$}%
    \setbox\tw@\hbox{\m@th$#1%
       \widehat{%
          \vrule\@width\z@\@height\ht\z@
          \vrule\@height\z@\@width\wd\z@}$}%
    \dp\tw@-\ht\z@
    \@tempdima\ht\z@ \advance\@tempdima2\ht\tw@ \divide\@tempdima\thr@@
    \setbox\tw@\hbox{%
       \raise\@tempdima\hbox{\scalebox{1}[-1]{\lower\@tempdima\box
\tw@}}}%
    {\ooalign{\box\tw@ \cr \box\z@}}}
\makeatother

\newtheoremstyle{xx}
  {4pt}
  {0pt}
  {\upshape}
  {\bfseries}
  {}
  { }
  {}
  
\makeatletter 
 \newtheoremstyle{myu}%
  {\upshape\item[ \indent\indent\bf\underline{\theorem@headerfont ##2:}]}%
\makeatother
\makeatletter 
 \newtheoremstyle{myn}%
  {\item[\hskip\labelsep \ \bf ##1 \theorem@headerfont ##2.]}%
\makeatother
\theoremstyle{myn}
\newtheorem{theoremn}{Theorem} 
\theoremstyle{myu}
{\upshape}
\newtheorem{x}[theoremn]{}

 \newtheoremstyle{mr}%
  {\upshape\item[ \indent{\theorem@headerfont ##2. \hspace{.2cm}}]}%
\makeatother
\theoremstyle{mr}
{\upshape}
\newtheorem{rf}[theoremn]{}

\title{\textbf{Analysis 101:\\
Curves and Length}}
\author{Garth Warner\\
Department of Mathematics\\
University of Washington}
\date{}	

\titleformat{\chapter}[display]
{\normalfont\filcenter\huge\bfseries}{}{0pt}{\large}

\titleformat{\chapter}[display]
{\normalfont\filcenter\huge\bfseries}{}{0pt}{\large}

\usepackage[OT2,OT1]{fontenc} 
\newcommand\cyr
{
\renewcommand\rmdefault{wncyr} 
\renewcommand\sfdefault{wncyss} 
\renewcommand\encodingdefault{OT2} 
\normalfont
\selectfont
}

\DeclareTextFontCommand{\textcyr}{\cyr}

\makeindex 
\begin{document}

\maketitle                              

\titlespacing*{\chapter}{0pt}{-50pt}{40pt}
\setlength{\parskip}{0.1em}
\pagenumbering{bychapter}
\setcounter{chapter}{0}
\pagenumbering{bychapter}

\begingroup
\fontsize{11pt}{11pt}\selectfont


\[
\textbf{ABSTRACT}
\]
\\

This is a systematic accounting of the classical theorems of Jordan and Tonelli,
as well as an introduction to the theory of the Weierstrass integral which in its definitive form is due to Cesari.
This is installment II of a four part discussion of certain aspects of Real Analysis: Functions of a Single Variable, Curves and Length, Functions of Several Variables, and Surfaces and Area.
\\[2cm]

\[
\textbf{ACKNOWLEDGEMENT}
\]

Many thanks to David Clark for his rendering the original transcript into AMS-LaTeX.  
Both of us also thank Judith Clare for her meticulous proofreading.
\newpage

\[
\textbf{CURVES AND LENGTH}
\]
\\


\hspace{2.1cm} \ \S1. \quad FUNDAMENTALS%
\\[-.26cm]

\hspace{2.1cm} \ \S2. \quad ESTIMATES%
\\[-.26cm]

\hspace{2.1cm} \ \S3. \quad EQUIVALENCES%
\\[-.26cm]

\hspace{2.1cm} \ \S4. \quad FR\'ECHET DISTANCE%
\\[-.26cm]

\hspace{2.1cm} \ \S5. \quad THE REPRESENTATION THEOREM %
\\[-.26cm]

\hspace{2.1cm} \ \S6. \quad INDUCED MEASURES %
\\[-.26cm]

\hspace{2.1cm} \ \S7. \quad TWO THEOREMS %
\\[-.26cm]

\hspace{2.1cm} \ \S8. \quad LINE INTEGRALS%
\\[-.26cm]

\hspace{2.1cm} \ \S9. \quad QUASI ADDITIVITY %
\\[-.26cm]

\hspace{2.1cm} \S10. \quad LINE INTEGRALS (bis) %
\\[-.26cm]

\hspace{2.1cm} \S11. \quad EXAMPLES %
\\[-.26cm]

\hspace{3.35cm}   REFERENCES
\\

\[
\]



\endgroup 

\chapter{
$\boldsymbol{\S}$\textbf{1}.\quad FUNDAMENTALS}
\setlength\parindent{2em}
\setcounter{theoremn}{0}
\renewcommand{\thepage}{\S1-\arabic{page}}

\begin{x}{\small\bf NOTATION} \ 
Given
\[
\un{x} = (x_1, \ldots, x_M) \in \R^M \quad (M = 1, 2, \ldots),
\]
put
\[
\norm{\un{x}} = (x_1^2 + \ldots + x_M^2)^\frac{1}{2},
\]
hence
\[
\abs{x_m} 
\ \leq \ 
\norm{\un{x}} 
\ \leq \ 
\abs{x_1} + \ldots + \abs{x_M} \quad (m = 1, \ldots, M).
\]
\\[-.25cm]
\end{x}

\begin{x}{\small\bf DEFINITION} \ 
A function $\un{f}:[a,b] \rightarrow \R^M$ is said to be a \un{curve} $C$, 
denoted $C \longleftrightarrow \un{f}$, where
\[
\un{f}(x) = (f_1(x), \ldots, f_M(x)) \quad (a \leq x \leq b).
\]
\\[-1.25cm]
\end{x}

\begin{x}{\small\bf EXAMPLE} \ 
Every function $f:[a,b] \rightarrow \R$ gives rise to a curve $C$ in $\R^2$, viz. the arrow 
\[
x \ra (x,f(x)).
\]
\\[-1.25cm]
\end{x}

\begin{x}{\small\bf DEFINITION} \ 
The \un{graph} of $C$, denoted $[C]$, is the range of $\un{f}$.
\\[-.25cm]
\end{x}

\begin{x}{\small\bf EXAMPLE} \ 
Take $M = 2$, let $k = 1, 2, \ldots$, and put
\[
\un{f_k}(x) 
\ = \  
\left( \sin^2 (kx),0\right) \qquad \left(0 \leq x \leq \frac{\pi}{2}\right).
\]
Then the $\un{f_k}$ all have the same range, i.e.,
\[
[C_1] 
\ = \  
[C_2] 
\ = \  
\ldots \qquad \text{ if } \quad C_k \longleftrightarrow \un{f_k}
\]
but the $C_k$ are different curves.
\\[-.25cm]
\end{x}

\begin{x}{\small\bf REMARK} \  
If $C$ is a continuous curve, then its graph $[C]$ is closed, 
bounded, connected, and uniformly locally connected.  
Owing to a theorem of Hahn 
and Mazurkiewicz, these properties are characteristic:  
Any such set is the graph of a continuous curve.  
So, e.g., a square in $\R^2$ is the graph of a continuous curve, a cube in $\R^3$ is the graph of a continuous curve etc.
\\[-.25cm]
\end{x}

\begin{x}{\small\bf DEFINITION} \  
The \un{length} of a curve $C$, denoted $\ell(C)$, is 
\[
\tT_{\un{f}}[a,b] 
\ \equiv \  
\sup\limits_{P \in \sP[a,b]} \ 
\sum\limits_{i=1}^n \ 
\norm{\un{f}(x_i) - \un{f}(x_{i-1})},
\]
$C$ being termed \un{rectifiable} if $\ell(C) < +\infty$.
\\[-.5cm]

[Note:  If $C$ is continuous and rectifiable, then $\forall \ \varepsilon > 0, \ \exists \ \delta > 0:$
\[
\norm{P} < \delta 
\implies 
\bigvee\limits_a^b \ (\un{f};P) 
\ \equiv \  
\sum\limits_{i=1}^n \ 
\norm{\un{f}(x_i) - \un{f}(x_{i-1})} 
\ > \  
\ell(C) - \varepsilon.]
\]
\\[-1cm]
\end{x}

\begin{x}{\small\bf LEMMA} \ 
Given a curve $C$, 
\[
\tT_{f_m}[a,b] 
\ \leq \ 
\ell(C) 
\ \leq \  
\tT_{f_1}[a,b] + \ldots +\tT_{f_M}[a,b] \qquad (1 \leq m \leq M).
\]
\\[-1cm]
\end{x}

\begin{x}{\small\bf SCHOLIUM} \ 
$C$ is rectifiable iff
\[
f_1 \in \BV[a,b], \ldots, f_M \in \BV[a,b].
\]
\\[-1cm]
\end{x}

\begin{x}{\small\bf THEOREM} \ 
Let
\begin{equation}
\begin{cases}
\ C_n \longleftrightarrow \un{f_n}:[a,b] \rightarrow \R^M \nonumber \\
\ C \hspace{0.15cm} \longleftrightarrow \un{f}:[a,b] \rightarrow \R^M \nonumber \\
\end{cases}
\end{equation}
and assume that $\un{f_n}$ converges pointwise to $\un{f}$ $-$then
\[
\ell(C) 
\ \leq \  
\liminf\limits_{n \rightarrow \infty} \ \ell(C_n).
\]
\\[-1cm]
\end{x}

A continuous curve
\[
\Gamma \longleftrightarrow \gamma:[a,b] \rightarrow \R^M
\]
is said to be a \un{polygonal line} 
(and \un{$\gamma$ quasi linear} in $[a,b]$) if there exists a $P \in \sP[a,b]$ in each segment of which $\un{\gamma}$ is linear or a constant.
\\[-.25cm]

\begin{x}{\small\bf DEFINITION} \ 
The \un{elementary length} $\ell_e(\Gamma)$ of $\Gamma$ is the sum of the lengths of these segments, hence $\ell_e(\Gamma) = \ell(\Gamma)$.
\\[-.25cm]
\end{x}

\begin{x}{\small\bf NOTATION} \ 
Given a continuous curve $C$, denote by $\Gamma(C)$ the set of all sequences
\[
\Gamma_n \longleftrightarrow \un{\gamma_n}:[a,b] \rightarrow \R^M
\]
of polygonal lines such that
\[
\gamma_n \rightarrow \un{f} \qquad (n \rightarrow \infty)
\]
uniformly in $[a,b]$.
\\[-.25cm]

Therefore
\[
\ell(C) 
\ \leq \  
\liminf\limits_{n \rightarrow \infty} \  \ell(\Gamma_n) 
\ = \  
\liminf\limits_{n \rightarrow \infty} \ \ell_e(\Gamma_n).
\]
On the other hand, by definition, there is some $\{\Gamma_n\} \in \Gamma(C)$ such that
\[
\ell_e(\Gamma_n) \rightarrow \ell(C) \qquad (n \rightarrow \infty).
\]
\\[-1cm]
\end{x}

\begin{x}{\small\bf SCHOLIUM} \ 
If $C$ is a continuous curve, then
\[
\ell(C) 
\ = \  
\inf\limits_{\{\Gamma_n\} \in \Gamma(C)} \ 
[\liminf\limits_{n \rightarrow \infty} \ell_e(\Gamma_n)].
\]
\\[-1cm]
\end{x}

\begin{x}{\small\bf REMARK} \ 
Let
\[
C \longleftrightarrow f:[a,b] \rightarrow \R^M.
\]
Assume:  $C$ is continuous and rectifiable $-$then $f$ can be decomposed as a sum
\[
f = f_{AC} + f_C
\]
where $f_{AC}$ is absolutely continuous and $f_C$ is continuous and singular.
\\[-.5cm]

\noindent
Therefore
\[
\ell(C) \ = \  \tT_{f_{AC}}[a,b] + \tT_{f_C}[a,b].
\]
\end{x}

\chapter{
$\boldsymbol{\S}$\textbf{2}.\quad ESTIMATES}
\setlength\parindent{2em}
\setcounter{theoremn}{0}
\renewcommand{\thepage}{\S2-\arabic{page}}

\begin{x}{\small\bf NOTATION} \ 
Write
\[
\tT_{\un{f}} [a,b]
\]
in place of $\ell(C)$.
\\[-.25cm]
\end{x}

\begin{x}{\small\bf DEFINITION} \ 
Assume that $C$ is rectifiable $-$then the \un{arc length function}
\[
s:[a,b] \ra \R
\]
 is defined by the prescription
\[
s(x) = \tT_{\un{f}}[a,x] \qquad (a \leq x \leq b).
\]

Obviously
\[
s(a) 
\ = \  
0, \quad s(b) 
\ = \  
\ell(C),
\]
and $s$ is an increasing function.
\\[-.25cm]
\end{x}

\begin{x}{\small\bf LEMMA} \ 
If $C$ is continuous and rectifiable, then $s$ is  continuous as are the $\tT_{f_m}[a,-] \quad (m = 1, \ldots, M)$.
\\[-.25cm]
\end{x}

\begin{x}{\small\bf LEMMA} \ 
If $C$ is continuous and rectifiable, then $s$ is absolutely continuous iff all the 
$\tT_{f_m}[a,-] \quad (m = 1, \ldots, M)$ are absolutely continuous, 
hence iff all the $f_m \quad (m = 1, \ldots, M)$ are absolutely continuous.
\\[-.25cm]
\end{x}

If $C$ is continuous and rectifiable, then the $f_m \in \BV[a,b]$, 
thus the derivatives $f_m^\prime$ exist almost everywhere in $[a,b]$ and are Lebesgue integrable.  
On the other hand, $s$ is an increasing function, 
thus it too is differentiable almost everywhere in $[a,b]$ and is Lebesgue integrable.
\\[-.25cm]


\begin{x}{\small\bf SUBLEMMA} \ 
The connection between $\un{f}^\prime$ and $s^\prime$ is given by the relation
\[
\norm{\un{f}^\prime}
\ \leq \   
s^\prime
\]
almost everywhere in $[a,b]$.
\\[-.5cm]

[For any subinterval $[\alpha,\beta] \subset [a,b]$,
\[
\norm{\un{f}(\beta) - \un{f}(\alpha)} 
\ \leq \   
s(\beta) - s(\alpha).]
\]
\\[-1.cm]
\end{x}

\begin{x}{\small\bf LEMMA} \ 
\[
\ell(C) 
\ = \  
s(b) - s(a) 
\ \geq \   
\int\limits_a^b s^\prime 
\ \geq \   
\int\limits_a^b \norm{\un{f}^\prime}.
\]
I.e.:  Under the assumption that $C$ is continuous and rectifiable,
\[
\ell(C) 
\ \geq \ 
\int\limits_a^b \ 
\norm{\un{f}^\prime}.
\]
\\[-1.cm]
\end{x}

\begin{x}{\small\bf THEOREM} \ 
\[
\ell(C) 
\ = \  
\int\limits_a^b \ 
\norm{\un{f}^\prime}
\]
iff all the $f_m$ $(m = 1, \ldots, M)$ are absolutely continuous.
\\[-.25cm]
\end{x}

This is established in the discussion to follow.
\\[-.25cm]

\qquad \textbullet \quad Suppose that the equality sign obtains, hence
\[
s(b) - s(a) 
\ = \  
\int\limits_a^b \ s^\prime.
\]
But also
\begin{align*}
&s(x) - s(a) \ \geq \  \int\limits_a^x \ s^\prime  
\\[15pt]
&s(b) - s(x) \ \geq \  \int\limits_x^b s^\prime.  
\end{align*}
If
\begin{align*}
&s(x) - s(a) \ > \  \int\limits_a^x s^\prime  
\\[15pt]
&s(b) - s(x) \ \geq \  \int\limits_x^b s^\prime, 
\end{align*}
then
\[
s(b) - s(a) \ > \  \int\limits_a^b s^\prime,
\]
a contradiction.  Therefore
\[
s(x) - s(a) \ = \  \int\limits_a^x s^\prime
\]

$
\implies \qquad s \in \AC[a,b]
$
\\[-.5cm]

$
\implies \qquad f_m \in \AC[a,b] \quad (m = 1, \ldots, M).
$
\\[-.25cm]

\qquad \textbullet \quad Consider the other direction, i.e., assume that the $f_m \in \AC[a,b]$, the claim being that
\[
\ell(C) = \int\limits_a^b \norm{f^\prime}.
\]
Given $P \in \sP[a,b]$, write

\allowdisplaybreaks
\begin{align*}
\sum\limits_{i=1}^n \ 
\norm{\un{f}(x_i) - \un{f}(x_{i-1})} \ 
&= \ 
\sum\limits_{i=1}^n \ 
\left[\sum\limits_{m=1}^M (f_m(x_i) - f_m(x_{i-1}))^2\right]^\frac{1}{2}    
\\[15pt]
&= \ 
\sum\limits_{i=1}^n \ 
\left[
\sum\limits_{m=1}^M \ 
\bigg( \hsx
\int\limits_{x_i}^{x_i} \ 
f_m^\prime
\bigg)^2
\right]^\frac{1}{2}  
\\[15pt]
&\leq \ 
\sum\limits_{i=1}^n  \ 
\int\limits_{x_i}^{x_i} \ 
\bigg( \hsx
\sum\limits_{m=1}^M (f_m^\prime)^2 
\bigg)^{1/2} 
\\[15pt]
&= \ 
\int\limits_a^b \ 
\norm{f^\prime}. 
\end{align*}
Taking the sup of the first term over all $P$ then gives
\[
\ell(C) 
\ \leq \  
\int\limits_a^b \ 
\norm{\un{f}^\prime} \quad (\leq \ell(C))
\]
$\implies$
\[
\ell(C) 
\ = \  
\int\limits_a^b \ 
\norm{\un{f}^\prime}.
\]
\\[-1cm]

\begin{x}{\small\bf \un{N.B.}} \ 
Under canonical assumptions,
\[
\bigg( \hsx
\bigg( \hsx
\int\limits_\X \phi_1
\bigg)^2 + \ldots + 
\bigg(
\int\limits_\X \ 
\phi_n
\bigg)^2 
\hsx\bigg)^\frac{1}{2} 
\ \leq \   
\int\limits_\X 
\left(
\phi_1^2 + \ldots + \phi_n^2
\right)^\frac{1}{2}.
\]
\\[-1cm]
\end{x}


\begin{x}{\small\bf RAPPEL} 
Suppose that $f \in \BV[a,b]$ $-$then for almost all $x \in [a,b]$,
\[
\abs{f^\prime(x)} 
\ = \  
\tT_{f^\prime}[a,x].
\]
\\[-1cm]
\end{x}

\begin{x}{\small\bf LEMMA} \ 
Suppose that $C$ is continuous and rectifiable $-$then
\[
s^\prime \ = \  \norm{\un{f}^\prime}
\]
almost everywhere in $[a,b]$.
\\[-.5cm]

PROOF \  
Since
\[
\norm{\un{f}^\prime} 
\ \leq \   
s^\prime,
\]
it suffices to show that
\[
s^\prime 
\ \leq \   
\norm{\un{f}^\prime}.
\]
Let $E_0 \subset [a,b]$ be the set of $x$ such that $\un{f}$ and $s$ are differentiable at $x$ 
and $s^\prime(x) > \norm{\un{f}^\prime(x)}$ and for $k = 1, 2, \ldots$, 
let $E_k$ be the set of $x \in E_0$ such that
\[
\frac{s(t_2) - s(t_1)}{t_2 - t_1} 
\ \geq \   
\frac{\norm{\un{f}(t_2) - \un{f}(t_1)}}{t_2 - t_1} + \frac{1}{k}
\]
for all intervals $[t_1,t_2]$ such that $x \in [t_1,t_2]$ and $0 < t_2 - t_1 \leq \ds\frac{1}{k}$.  
So, by construction,
\[
E_0 
\ = \  
\bigcup\limits_{k=1}^\infty E_k
\]
and matters reduce to establishing that $\forall \ k, \lambda(E_k) = 0$.  To this end, let $\varepsilon > 0$ and choose $P \in \sP[a,b]:$
\[
\sum\limits_{i=1}^n 
\norm{\un{f}(x_i) - \un{f}(x_{i-1})} 
\ > \  
\tT_{\un{f}} [a,b] - \varepsilon.
\]
Expanding $P$ if necessary, it an be assumed without loss of generality that
\[
0 
\ < \  
x_i - x_{i-1} 
\ \leq \   
\frac{1}{k} \qquad (i = 1, \ldots, n).
\]
For each i, either $[x_{i-1},x_i] \cap E_k \neq \emptyset$ and then
\[
s(x_i) - s(x_{i-1}) 
\ \geq \   
\norm{\un{f}(x_i) - \un{f}(x_{i-1})} + \frac{x_i - x_{i-1}}{k},
\]
or $[x_{i-1}, x_i] \cap E_k = \emptyset$ and then
\[
s(x_i) - s(x_{i-1}) 
\ \geq \  
\norm{\un{f}(x_i) - \un{f}(x_{i-1})}.
\]
Consequently

\allowdisplaybreaks
\begin{align*}
\tT_{\un{f}} [a,b] \ 
&= \ 
s(b)  
\\[15pt]
&=\ 
s(x_n)  
\\[15pt]
&= \ 
\sum\limits_{i=1}^n (s(x_i) - s(x_{i-1})) \qquad 
(s(x_0) = s(a) = 0 ) 
\\[15pt]
&\geq \ 
\sum\limits_{i=1}^n \ 
\norm{\un{f}(x_i) - \un{f}(x_{i-1})} + \frac{1}{k} \lambda^*(E_k)  
\\[15pt]
&\geq \ 
\tT_{\un{f}} [a,b] - \varepsilon + \frac{1}{k} \lambda^*(E_k)  
\end{align*}

$\implies$
\[
\lambda^*(E_k) 
\ \leq \   
k \varepsilon 
\quad \implies \quad 
\lambda(E_k) 
\ = \  
0 
\quad (\varepsilon \downarrow 0).
\]
\\[-1cm]
\end{x}

\begin{x}{\small\bf THEOREM} \ 
Suppose that $C$ is continuous and rectifiable.  
Assume:  $M > 1$ $-$then the $M$-dimensional Lebesgue measure of $[C]$ is equal to 0.
\\[-.25cm]
\end{x}

\begin{x}{\small\bf NOTATION} \ 
Let
\[
C \longleftrightarrow \un{f}:[a,b] \ra \R^M
\]
be a continuous curve.  
Given $\un{x} \in [C]$, let $N(\un{f};\un{x})$ be the number of points $x \in [a,b]$ (finite or infinite) 
such that $f(x) = \un{x}$ and let $N(\un{f};-) = 0$ in the complement $\R^M - [C]$ of $[C]$.
\\[-.25cm]
\end{x}

\begin{x}{\small\bf THEOREM} \ 
\[
\ell(C) 
\ = \  
\int\limits_{\R^M} \  N(\un{f};-) \ \td H^1.
\]

[Note:\ 
$H^1$ is the $1$-dimensional Hausdorff outer measure in $\R^M$ and
\[
H^1([C]) 
\ = \  
\int\limits_{\R^M} \chi_{[C]}\  \td H^1 
\ \leq \  
\int\limits_{\R^M} N(\un{f};-) \ \td H^1,
\]
i.e.,
\[
H^1([C]) 
\ \leq \  
\ell(C)
\]
and it can happen that
\[
H^1([C]) 
\ < \  
\ell(C).]
\]
\\[-1cm]
\end{x}

\begin{x}{\small\bf \un{N.B.}} \ 
If $\un{f}$ is one-to-one, then
\[
N(\un{f};-) 
\ = \  
\chi_{[C]}
\]
and when this is so,
\[
H^1([C]) = \ell(C).
\]
\end{x}

\chapter{
$\boldsymbol{\S}$\textbf{3}.\quad EQUIVALENCES}
\setlength\parindent{2em}
\setcounter{theoremn}{0}
\renewcommand{\thepage}{\S3-\arabic{page}}

\qquad 
In what follows, by \un{interval} we shall understand a finite closed interval $\subset \R$.
\\[-.5cm]

[Note: If $I, J$ are intervals and if $\partial I = \{a,b\}$, $\partial J = \{c,d\}$, 
then the agreement is that a homeomorphism $\phi:I \ra J$ is sense preserving, i.e., sends $a$ to $c$ and $b$ to $d$.]
\\[-.25cm]

\begin{x}{\small\bf DEFINITION} \ 
Suppose given intervals $I$, $J$ and curves 
$\un{f}: I \ra \R^M$, $\un{g}: J \ra \R^M$ $-$then $\un{f}$ 
and 
$\un{g}$ are said to be \un{Lebesgue equivalent} if there exists a homeomorphism 
$\phi:I \ra J$ such that $\un{f} = \un{g} \circ \phi$.
\\[-.25cm]
\end{x}

\begin{x}{\small\bf LEMMA} \ 
If
\[
\begin{cases}
\ \un{f}: [a,b] \ra \R^M  
\\[4pt]
\ \un{g}: [c,d] \ra \R^M  
\end{cases}
\]
are Lebesgue equivalent and if
\[
\begin{cases}
\ C \longleftrightarrow \un{f}  
\\[4pt]
\ D \longleftrightarrow \un{g}
\end{cases}
,
\]
then
\[
\ell(C) 
\ = \  
\ell(D).
\]

PROOF \ 
The homeomorphism $\phi:[a,b] \ra [c,d]$ induces a bijection
\[
\begin{cases}
\ \sP[a,b] &\ra \sP[c,d]  \\
\ \ P &\ra Q
\end{cases}
.
\]
Therefore
\allowdisplaybreaks
\begin{align*}
\ell(C) 
&= \sup\limits_{P \in \sP[a,b]} \ 
\sum\limits_{i=1}^n \ 
\norm{\un{f}(x_i) - \un{f}(x_{i-1})}  
\\[15pt]
&= \sup\limits_{P \in \sP[a,b]} \ 
\sum\limits_{i=1}^n \ 
\norm{\un{g}(\phi(x_i)) - \un{g}(\phi(x_{i-1}))}  
\\[15pt]
&= \sup\limits_{Q \in \sP[c,d]} \ 
\sum\limits_{i=1}^n \ 
\norm{\un{g}(y_i) - \un{g}(y_{i-1})}  
\\[15pt]
&= \ell(D).  \\
\end{align*}
\\[-1.cm]
\end{x}

\begin{x}{\small\bf DEFINITION} \ 
Suppose given intervals $I$, $J$ and curves 
$\un{f}: I \ra \R^M$, $\un{g}: J \ra \R^M$ $-$then $\un{f}$ and $\un{g}$ are said to be 
\un{Fr\'echet equivalent} if for every $\varepsilon > 0$ there exists a homeomorphism $\phi:I \ra J$ such that
\[
\norm{\un{f}(x) - \un{g}(\phi(x))} 
\ < \  
\varepsilon \qquad (x \in I).
\]
\\[-1.cm]
\end{x}

\begin{x}{\small\bf REMARK} \ 
It is clear that two Lebesgue equivalent curves are Fr\'echet equivalent 
but two Fr\'echet equivalent curves need not be Lebesgue equivalent.
\\[-.25cm]
\end{x}

\begin{x}{\small\bf LEMMA} \ 
If
\[
\begin{cases}
\ \un{f}: [a,b] \ra \R^M  
\\[4pt]
\ \un{g}: [c,d] \ra \R^M  
\end{cases}
\]
are Fr\'echet equivalent and if
\[
\begin{cases}
\ C \longleftrightarrow \un{f}  
\\[4pt]
\ D \longleftrightarrow \un{g}
\end{cases}
,
\]
then
\[
\ell(C) \ = \  \ell(D).
\]

PROOF \ 
For each $n = 1, 2, \ldots$, there is a homeomorphism $\phi_n[a,b] \ra [c,d]$ such that $\forall \ x \in [a,b]$, 
\[
\norm{\un{f}(x) - \un{g}(\phi_n(x))} 
\ < \  
\frac{1}{n}.
\]
Put $\un{f}_{\hsx n} = \un{g} \circ \phi_{n}$, hence $\un{f}_{\hsx n}$ is Lebesgue equivalent to $\un{g}$ 
(viz. $\un{g} \circ\phi_n = \un{g} \circ \phi_n \ldots$ ), 
thus if
\[
C_n \longleftrightarrow \un{f}_{\hsx n}, \quad D \longleftrightarrow \un{g}, 
\]
then from the above
\[
\ell(C_n) 
\ = \  
\ell(D).
\]
But $\forall \ x \in [a,b]$,
\[
\norm{\un{f}(x) - \un{f}_{\hsx n}(x)} 
\ < \  
\frac{1}{n},
\]
i.e., $\un{f}_{\hsx n} \ra \un{f}$ pointwise, so
\begin{align*}
\ell(C) \ 
&\leq \ 
\liminf\limits_{n \ra \infty} \ \ell(C_n) 
\\[11pt]
&\leq \ 
\liminf\limits_{n \ra \infty} \ \ell(D) 
\\[11pt]
&\leq \ 
\ell(D). 
\end{align*}
Analogously
\[
\ell(D) 
\ \leq \ 
\ell(C).
\]
Therefore
\[
\ell(C) 
\ = \  
\ell(D).
\]
\end{x}

\chapter{
$\boldsymbol{\S}$\textbf{4}.\quad FR\'ECHET DISTANCE}
\setlength\parindent{2em}
\setcounter{theoremn}{0}
\renewcommand{\thepage}{\S4-\arabic{page}}

\qquad 
Let
\[
\begin{cases}
\ C \longleftrightarrow \un{f}: [a,b] \ra \R^M  
\\[8pt]
\ D \longleftrightarrow \un{g}: [c,d] \ra \R^M  
\end{cases}
\]
be two continuous curves.
\\[-.25cm]

\begin{x}{\small\bf NOTATION} \  
$\sH$ is the set of all homeomorphisms
\[
\phi: [a,b] \ra [c,d] \qquad (\phi(a) = c, \ \phi(b) = d).
\]
Given $\phi \in \sH$ , the expression 
\[
\norm{\un{f}(x) - \un{g}(\phi(x))} \qquad (a \leq x \leq b)
\]
has an absolute maximum $M(\un{f},\un{g};\phi)$.
\\[-.25cm]
\end{x}

\begin{x}{\small\bf DEFINITION} \ 
The \un{Fr\'echet distance} between $C$ and $D$, denoted $\norm{C,D}$, is 
\[
\inf\limits_{\phi \in \sH} M(\un{f},\un{g});\phi).
\]

[Note: \ 
In other words, $\norm{C,D}$ is the infimum of all numbers $\varepsilon \geq 0$ with the property that there exists a homeomorphism $\phi \in \sH$ such that
\[
\norm{\un{f}(x) - \un{g}(\phi(x))}
\ \leq \   
\varepsilon
\]
for all $x \in [a,b].]$
\\[-.25cm]
\end{x}

\begin{x}{\small\bf \un{N.B.}} \ 
If $\norm{C,D} < \varepsilon$, then there exists a $\phi \in \sH$ such that
\[
M(\un{f},\un{g};\phi) 
\ < \  
\varepsilon.
\]
\\[-1cm]
\end{x}

\begin{x}{\small\bf LEMMA} \ 
Let $C$, $D$, $C_0$ be continuous curves $-$then
\\[-.25cm]

\qquad\qquad (i) \quad 
$\norm{C,D} \ \geq \ 0;$
\\[-.25cm]

\qquad\qquad (ii) \quad 
$\norm{C,D} \ = \  \norm{D,C}; $ 
\\[-.25cm]

\qquad\qquad (iii) \quad 
$\norm{C,D} \ \leq \ \norm{C,C_0} + \norm{C_0,D};$ 
\\[-.25cm]

\qquad\qquad (iv) \quad 
$\norm{C,D} \ = \  0 \ \text{iff $C$ and $D$ are Fr\'echet equivalent.}$  
\\[-.25cm]

Therefore Fr\'echet distance is a premetric on the set of continuous curves with values in $\R^M$.
\\[-.25cm]
\end{x}

\begin{x}{\small\bf THEOREM} \ 
Let
\[
\begin{cases}
\ C_n \longleftrightarrow \un{f_n}: [a_n,b_n] \ra \R^M \qquad (n = 1, 2, \ldots) 
\\[11pt]
\ C \longleftrightarrow \un{f}: [a,b] \ra \R^M  
\end{cases}
\]
be continuous curves.  
Assume:
\[
\norm{C_n,C} \ra 0 \qquad (n \ra \infty).
\]
Then
\[
\qquad \ell(C) 
\ \leq \  
\liminf\limits_{n \ra \infty} \ \ell(C_n).
\]
\\[-1cm]

PROOF \ 
For every $n$, there is a homeomorphism
\[
\phi_n:[a,b] \ra [a_n,b_n] \qquad (\phi_n(a) = a_n, \ \phi_n(b) = b_n)
\]
such that for all $x \in [a,b]$,
\[
\norm{\un{f}(x) - \un{f_n}(\phi_n(x))} 
\ < \  
\norm{C,C_n} + \frac{1}{n}.
\]
Let
\[
D_n \longleftrightarrow \un{f_n} \circ \phi_n: [a,b] \ra \R^M.
\]
Then pointwise
\[
\un{f_n} \circ \phi_n \ra \un{f}
\]

$\implies$
\[
\ell(C) 
\ \leq \   
\liminf\limits_{n \ra \infty} \ell(D_n).
\]
But \quad $\ell(D_n) = \ell(C_n)$, hence
\[
\ell(C) 
\ \leq \   
\liminf\limits_{n \ra \infty} \ell(C_n).
\]
\\[-1.5cm]
\end{x}

\qquad In the set of continuous curves, 
introduce an equivalence relation by stipulating that 
$C$ and $D$ are equivalent provided $C$ and $D$ are Fr\'echet equivalent.  
The resulting set $\sE_F$ of equivalence classes is then a metric space: If
\[
\begin{cases}
\ \{C\} \in \sE_F  
\\[8pt]
\ \{D\} \in \sE_{F}
\end{cases}
,
\]
then
\[
\norm{\{C\},\{D\}} 
\ = \  
\norm{C,D}.
\]
\\[-1cm]

\begin{x}{\small\bf \un{N.B.}} \ %
If $C$, $C^\prime$ are Fr\'echet equivalent and if $D$, $D^\prime$ are Fr\'echet equivalent, then
\begin{align*}
\norm{C,D} \ 
&\leq \ 
\norm{C,C^\prime} + \norm{C^\prime,D} 
\\[8pt]
&\leq \ 
\norm{C^\prime,D} 
\\[8pt]
&\leq \ 
\norm{C^\prime,D^\prime} + \norm{D^\prime,D} 
\\[8pt]
&= \ 
\norm{C^\prime,D^\prime}
\end{align*}
and in reverse
\[
\norm{C^\prime,D^\prime}
\ \leq \  
\norm{C,D}.
\]
So
\[
\norm{C,D} 
\ = \  
\norm{C^\prime,D^\prime}.
\]
\end{x}

\chapter{
$\boldsymbol{\S}$\textbf{5}.\quad THE REPRESENTATION THEOREM}
\setlength\parindent{2em}
\setcounter{theoremn}{0}
\renewcommand{\thepage}{\S5-\arabic{page}}

\qquad  
Assume:
\[
C \longleftrightarrow \un{f}: [a,b] \ra \R^M
\]
is a curve which is continuous and rectifiable.
\\[-.25cm]

\begin{x}{\small\bf THEOREM} \ %
There exists a continuous curve
\[
D \longleftrightarrow \un{g}: [c,d] \ra \R^M
\]
with the property that
\[
\ell(D) 
\ = \  
\ell(C) \qquad (< +\infty)
\]
and
\[
\ell(D) 
\ = \  
\int\limits_c^d \ \norm{g^\prime},
\]
where $g_1, \ldots, g_M$ are absolutely continuous and in addition $\un{f}$ and $\un{g}$ are Fr\'echet equivalent.
\end{x}

Take $\ell(C) > 0$ and define $\un{g}$ via the following procedure.  
In the first place, the domain $[c,d]$ of $\un{g}$ is going to be the interval $[0,\ell(C)]$.  
This said, note that $s(x)$ is contant in an interval $[\alpha,\beta]$ iff $\un{f} (x)$ is constant there as well.  
Next, for each point $s_0$ $(0 \leq s_0 \leq \ell(C))$ there is a maximal interval 
$\alpha \leq x \leq \beta \quad (a \leq \alpha \leq \beta \leq b)$ with $s(x) = s_0$. \ 
Definition: \ 
$\un{g}(s_0) = \un{f}(x) \quad (\alpha \leq x \leq \beta)$.
\\

\begin{x}{\small\bf LEMMA} \ 
\[
\begin{cases}
\ \un{g}(s_{0^-}) = \un{g}(s_0) \qquad (0 < s_0 \leq \ell(C))  
\\[4pt]
\ \un{g}(s_{0^+}) = \un{g}(s_0) \qquad (0 \leq s_0 < \ell(C))
\end{cases}
.
\]

Therefore
\[
\un{g}: [c,d] \ra \R^M
\]

is a continuous curve.
\\[-.25cm]
\end{x}

\begin{x}{\small\bf SUBLEMMA} \ 
Suppose that $\phi_n: [A,B] \ra [C,D]$ $(n = 1, 2, \ldots)$ converges uniformly to 
$\phi: [A,B] \ra [C,D]$.  Let $\Phi:[C,D] \ra \R^M$ be a continuous function 
$-$then $\Phi \circ \phi_n$ converges uniformly to $\Phi \circ \phi$.
\\[-.25cm]

PROOF \  
Since $\Phi$ is uniformly continuous, given $\varepsilon > 0$, $\exists \ \delta > 0$ such that
\[
\abs{u - v} 
\ < \  
\delta \ \implies \ \norm{\Phi(u) - \Phi(v)} 
\ < \  
\varepsilon \qquad (u,v \in [C,D]).
\]
Choose N:
\[
n 
\ \geq \  
N 
\ \implies \ 
\abs{\phi_n(x) - \phi(x)} 
\ < \  
\delta \qquad (x \in [A,B]).
\]
Then
\[
\norm{\Phi(\phi_n(x)) - \Phi(\phi(x))} 
\ < \  
\varepsilon .
\]
\\[-1cm]
\end{x}

\begin{x}{\small\bf LEMMA} \ 
$\un{f}$ and $\un{g}$ are Fr\'echet equivalent.
\\[-.25cm]

PROOF \  
Approximate $s$ by quasilinear, strictly increasing functions $s_n(x)$ 
$(a \leq x \leq b)$ with $s_n(a) = 0$, $s_n(b) = \ell(C)$ and
\[
\abs{s_n(x) - s(x)} 
\ < \  
\frac{1}{n} \qquad (n = 1, 2, \ldots).
\]
Then
\[
s_n: [a,b] \ra [0,\ell(C)]
\]
converges uniformly to
\[
s:[a,b] \ra [0,\ell(C)]
\]
and
\[
\un{g}: [0,\ell(C)] \ra \R^M
\]
is continuous, so
\[
\un{g} \circ s_n \ra \un{g} \circ s
\]
uniformly in $[a,b]$, thus $\forall \ \varepsilon > 0$, $\exists \ N: n \geq N$
\[
\implies \quad \norm{\un{g}(s_n(x)) - \un{g}(s(x))} 
\ < \  
\varepsilon \qquad (a \leq x \leq b)
\]
or still, 
\[
 \norm{\un{f}(x) - \un{g}(s_n(x))} 
 \ < \  
 \varepsilon \qquad (a \leq x \leq b).
\]
Since the $s_n$ are homeomorphisms, it follows that $\un{f}$ and $\un{g}$ are Fr\'echet equivalent.
\\[-.25cm]
\end{x}

\begin{x}{\small\bf LEMMA} \ 
\[
0 \leq u 
\ < \  
v 
\ \leq \  
\ell(C)
\]

\qquad
$\implies$
\[
\norm{\un{g}(v) - \un{g}(u)} 
\ = \   
v - u
\]

\qquad
$\implies$
\[
\abs{g_m(v) - g_m(u)} 
\ \leq \  
v - u \qquad (1 \leq m \leq M).
\]
Consequently $g_1, \ldots, g_M$ are absolutely continuous (in fact, Lipschitz).
\\[-.25cm]
\end{x}

\begin{x}{\small\bf LEMMA} \ 
\[
\ell(C) 
\ = \  
\ell(D) 
\ = \  
\int\limits_0^{\ell(D)} \  \norm{\un{g}^\prime},
\]
where $\norm{\un{g}^\prime} \leq 1$.
\\

\qquad So
\begin{align*}
0 \ 
&= \ 
\ell(D) - \int\limits_0^{\ell(D)} \ 
\norm{\un{g}^\prime} 
\\[15pt]
&= \ 
\int\limits_0^{\ell(D)} \ 
1 - \int\limits_0^{\ell(D)} \ 
\norm{\un{g}^\prime} 
\\[15pt]
&= \ 
\int\limits_0^{\ell(D)} \ 
(1 -  \norm{\un{g}^\prime} )
\end{align*}
implying thereby that $\norm{\un{g}^\prime} = 1$ almost everywhere.
\end{x}

\chapter{
$\boldsymbol{\S}$\textbf{6}.\quad INDUCED MEASURES}
\setlength\parindent{2em}
\setcounter{theoremn}{0}
\renewcommand{\thepage}{\S6-\arabic{page}}

\begin{x}{\small\bf NOTATION} \ 
$\tBO[a,b]$ is the set of Borel subsets of $[a,b]$.
\\[-.25cm]
\end{x}

Let
\[
C \longleftrightarrow \un{f}:[a,b] \ra \R^M
\]
be a curve, continuous and rectifiable.
\\[-.25cm]

\begin{x}{\small\bf LEMMA} \ 
The interval function defined by the rule
\[
[c,d] \ra s(d) - s(c) \qquad ([c,d] \subset [a,b])
\]
can be extended to a measure $\mu_C$ on $\tBO[a,b]$.
\\[-.25cm]
\end{x}

\begin{x}{\small\bf LEMMA} \ 
For $m = 1, \ldots, M$, the interval function defined by the rule
\[
[c,d] \ra T_{f_m}[c,d] \qquad ([c,d] \subset [a,b])
\]
can be extended to a measure $\mu_m$ on $\tBO[a,b]$.
\\[-.25cm]
\end{x}

\begin{x}{\small\bf FACT} \ 
Given $S \in \tBO[a,b]$, 
\[
\mu_m(S) \leq \mu_C(S) \leq \mu_1(S) + \ldots + \mu_M(S).
\]
\\[-1.5cm]
\end{x}

\begin{x}{\small\bf LEMMA} \ 
For $m = 1, \ldots, M$, the interval functions defined by the rule
\[
\begin{cases}
\ [c,d] \ra T_{f_m}^+[c,d]  
\\[8pt]
\ [c,d] \ra T_{f_m}^-[c,d]   
\end{cases}
\qquad ([c,d] \subset [a,b])
\]
can be extended to measures
\quad
$
\begin{cases}
\ \mu_m^+  
\\[4pt]
\ \mu_m^-  
\end{cases}
$
on $\tBO[a,b]$.
\\[-.25cm]
\end{x}


\begin{x}{\small\bf NOTATION} \ 
Put
\[
\nu_m 
\ = \  
\mu_m^+ - \mu_m^- \qquad (m = 1, \ldots, M).
\]
[Thus $\nu_m$ is a countably additive, totally finite set function on $\tBO[a,b].]$
\\[-.25cm]
\end{x}

\begin{x}{\small\bf RECOVERY PRINCIPLE} \ 
For any $S \in \tBO[a,b]$, 
\[
\mu_C(S) 
\ = \  
\sup\limits_{\{P\}} \ 
\sum\limits_{E \in P} \ 
\bigg\{ 
\sum\limits_{m=1}^M \ 
\nu_m(E)^2
\bigg\}^\frac{1}{2},
\]
where the supremum is taken over all partitions $P$ of $S$ into disjoint Borel measurable sets $E$.
\\[-.25cm]
\end{x}

\begin{x}{\small\bf FACT} \ 
The set functions $\mu_m, \  \mu_m^+,  \ \mu_m^-,  \ \nu_m$ are absolutely continuous w.r.t. $\mu_C$.
\\[-.25cm]
\end{x}

\begin{x}{\small\bf NOTATION} \ 
The corresponding Radon-Nikodym derivatives are denoted by
\[
\beta_m = \frac{\td\mu_m}{\td\mu_C}, \quad
\begin{cases}
\ \ds\beta^+_m = \frac{\td\mu_m^+}{\td\mu_C}  
\\[15pt]
\ \ds\beta^-_m = \frac{\td\mu_m^-}{\td\mu_C}  
\end{cases}
, \ \ds\theta_m = \frac{\td\nu_m}{\td\mu_C}.
\]
\\[-1cm]
\end{x}

\begin{x}{\small\bf CONVENTION} \ 
The term almost everywhere (or measure 0) will refer to the measure space
\[
([a,b], \tBO[a,b], \mu_C).
\]
\\[-1cm]
\end{x}

\begin{x}{\small\bf FACT} \ 
\[
\begin{cases}
\ \beta_m = \beta_m^+ + \beta_m^-  
\\[8pt]
\ \theta_m = \beta_m^+ - \beta_m^-  
\end{cases}
\qquad (m = 1, \ldots, M)
\]
almost everywhere.
\\[-.25cm]
\end{x}

\begin{x}{\small\bf NOTATION} \ 
Let
\[
\un{\theta} = (\theta_1, \ldots, \theta_M).
\]
[Note:  By definition,
\[
\norm{\un{\theta}(x)} = (\theta_1(x)^2 + \cdots+ \theta_M(x)^2)^\frac{1}{2}.]
\]
\\[-1cm]
\end{x}

\begin{x}{\small\bf NOTATION} \ 
Given a linear orthogonal transformation $\lambda: \R^M \ra \R^M$, let $\overline{C} = \lambda \hsy C$.
\\[-.25cm]
\end{x}

\begin{x}{\small\bf \un{N.B.}} \ 
\[
\mu_{\overline{C}} = \mu_C.
\]
\\[-1cm]
\end{x}

\begin{x}{\small\bf LEMMA} \ 
\[
(\overline{\nu}_1, \ldots, \overline{\nu}_M) = \lambda (\nu_1, \ldots, \nu_M).
\]
\\[-1cm]
\end{x}

\begin{x}{\small\bf APPLICATION} \ 
\[
(\overline{\theta}_1, \ldots, \overline{\theta}_M) = \lambda (\theta_1, \ldots, \theta_M)
\]
almost everywhere.
\\[-.5cm]

[Differentiate the preceding relation w.r.t. $\mu_{\overline{C}} = \mu_C.]$
\\[-.25cm]
\end{x}

\begin{x}{\small\bf LEMMA} \ 
\[
\abs{\theta_m} 
\ \leq \ 
1 \qquad  (m = 1, \ldots, M)
\]
almost everywhere, so
\[
\norm{\un{\theta}} 
\ \leq \  
M^\frac{1}{2}
\]
almost everywhere.
\\[-.25cm]
\end{x}

\begin{x}{\small\bf THEOREM} \ 
\[
\norm{\un{\theta}} 
\ = \  
1
\]
almost everywhere.
\\[-.25cm]

PROOF \ 
Let $0 < \delta < 1$ and let
\[
S 
\ = \ 
\{x: \norm{\un{\theta}(x)} < 1 - \delta\}.
\]
Then
\[
\mu_C(S) 
\ = \  
\sup\limits_{\{P\}} \ 
\sum\limits_{E \in P} \ 
\bigg\{ 
\sum\limits_{m=1}^M \ 
\nu_m(E)^2
\bigg\}^\frac{1}{2}.
\]
But
\begin{align*}
\nu_m(E) \ 
&= \ 
\int\limits_E \ 
\frac{\td \nu_m}{\td\mu_C} \ \td\mu_C 
\\[15pt]
&= \ 
\int\limits_E \ 
\theta_m \ \td\mu_C.
\end{align*}
Therefore

\allowdisplaybreaks
\begin{align*}
\bigg\{ 
\sum\limits_{m=1}^M \ \nu_m(E)^2 
\bigg\}^\frac{1}{2} \ 
&= \ 
\bigg\{ 
\sum\limits_{m=1}^M \ 
\bigg(
\int\limits_E \ 
\theta_m  \ \td\mu_C
\bigg)^2 
\bigg\}^\frac{1}{2} 
\\[15pt]
&\leq \ 
\int\limits_E \ 
\bigg\{  
\sum\limits_{m=1}^M  \ \theta_m^2 
\bigg\}^\frac{1}{2} 
\ \td\mu_C 
\\[15pt]
&= \ 
\int\limits_E \norm{\un{\theta}(x)} 
\ \td\mu_C 
\\[15pt]
&\leq \ 
(1 - \delta) \int\limits_E 
\ \td\mu_C
\\[15pt]
&= \ 
(1 - \delta) \hsy \mu_C(E).
\end{align*}
Since
\[
S = \coprod E,
\]
it follows that
\[
\sum\limits_{E \in P} \ 
\bigg\{ 
\sum\limits_{m=1}^M \ 
\nu_m(E)^2
\bigg\}^\frac{1}{2}
\ \leq \ 
 (1 - \delta) \hsx \mu_C(S).
\]
Taking the supremum over the $P$ then implies that
\[
\mu_C(S) 
\ \leq \  
(1 - \delta) \hsx \mu_C(S),
\]
thus $\mu_C(S) = 0$ and $\norm{\un{\theta}(x)} \geq 1$ almost everywhere 
(let $\ds\delta = \frac{1}{2}, \frac{1}{3}, \ldots$).  
To derive a contradiction, take $M \geq 2$ and suppose that 
$\norm{\un{\theta}(x)} \geq 1 + \delta > 1$ 
on some set $T$ such that $\mu_C(T) > 0$ 
$-$then for some vector
\[
\un{\xi} 
\ = \  
(\xi_1, \ldots, \xi_M) \in \R^M \qquad \left(\norm{\un{\xi}} = 1\right),
\]
the set
\[
T(\un{\xi}) 
\ = \  
\{x \in T: \normxx{\frac{\un{\theta}(x)}{\norm{\un{\theta}(x)}} - \un{\xi}} < \frac{\delta}{M^2}
\]
has measure $\mu_C (T(\un{\xi})) > 0$ $($see below$)$.  
Let
\[
\un{\lambda_j} 
\ = \  
 (\lambda_{j_1}, \ldots, \lambda_{j_M}) \qquad (j = 2, \ldots, M)
\]
be unit vectors such that
\\[-.25cm]
\end{x}

\[
\lambda = 
\begin{bmatrix}
\xi_1, &\ldots, &\xi_M\\
\lambda_{21}, &\ldots, &\lambda_{2M}\\
\vdots &\vdots &\vdots \\
\lambda_{M1}, &\ldots, &\lambda_{MM}\\
\end{bmatrix}
\]
is an orthogonal matrix.  
Viewing $\lambda$ as a linear orthogonal transformation, form as above $\overline{C} = \lambda C$, hence
\[
(\overline{\theta}_1, \ldots, \overline{\theta}_M) = \lambda(\theta_1, \ldots, \theta_M).
\]
On $T(\un{\xi})$,
\begin{align*}
\abs{\overline{\theta}_j} \ 
&= \ 
\abs{\lambda_{j_1} \theta_1 + \cdots + \lambda_{j_M} \theta_M}
\\[11pt]
&\leq \ 
\norm{\un{\theta}} \frac{\delta}{M^2}
\\[11pt]
&\leq \ 
M^\frac{1}{2} \frac{\delta}{M^2} 
\\[11pt]
&\leq \ 
M\frac{\delta}{M^2}
\\[11pt]
&= \ 
\frac{\delta}{M},
\end{align*}
while
\[
\norm{\overline{\un{\theta}}} 
\ \leq \  
\abs{\overline{\theta}_1} + \ldots + \abs{\overline{\theta}_M}
\]

\qquad
$\implies$
\begin{align*}
\abs{\overline{\theta}_1} \ 
&\geq \norm{\overline{\un{\theta}}} - \abs{\overline{\theta}_2} - \cdots - \abs{\overline{\theta}_M}
\\[11pt]
 &\geq \ 
 (1 + \delta) - (M - 1) \frac{\delta}{M}
\\[11pt]
 &= \ 
 1 + \frac{\delta}{M}.
\end{align*}
However
\[
\abs{\overline{\theta}_1} 
\ \leq \  
1,
\]
so we have a contradiction.
\\[-.25cm]

\begin{x}{\small\bf \un{N.B.}} \ 
Let $\{\un{\xi_n}: n \in \N\}$ be a dense subset of the unit sphere in $\tU(M)$ in $\R^M$ $($thus $\forall \ n$, $\norm{\un{\xi_n}} = 1)$.  
Given a point $x \in T$, pass to
\[
\frac{\un{\theta}(x)}{\norm{\un{\theta}(x)}} \in \tU(M).
\]
Then there exists a $\un{\xi_{n_x}}:$
\[
\norm{\frac{\un{\theta}(x)}{\norm{\un{\theta}(x)}} - \un{\xi_{n_x}}} 
\ < \  
\frac{\delta}{M^2},
\]

\noindent
a point in the $\ds\frac{\delta}{M^2}$ $-$ neighborhood of 
\[
\frac{\un{\theta}(x)}{\norm{\un{\theta}(x)}}
\]
in $\tU(M)$.  Therefore
\[
T 
\ = \  
\bigcup\limits_{n=1}^\infty \ T(\un{\xi_n})
\]

\qquad
$\implies$
\[
0 \ < \  
\mu_C(T) 
\ \leq \ 
\sum\limits_{n=1}^\infty \ \mu_C(T(\un{\xi_n}))
\]

\qquad
$\implies$ $\exists \ n:$
\[
\mu_C(T(\un{\xi_n})) \ > \  0.
\]
\end{x}

\chapter{
$\boldsymbol{\S}$\textbf{7}.\quad TWO THEOREMS}
\setlength\parindent{2em}
\setcounter{theoremn}{0}
\renewcommand{\thepage}{\S7-\arabic{page}}

\qquad 
Let
\[
C \longleftrightarrow \un{f}: [a,b] \rightarrow \R^M
\]
be a curve, continuous and rectifiable.
\\[-.25cm]

Let $P \in \sP[a,b]$, say
\[
P: a = x_0 < x_1 < \cdots < x_n = b.
\]
\\[-1cm]

\begin{x}{\small\bf DEFINITION} \ 
Let $i = 1, \ldots, n$ and for $m = 1, \ldots, M$ let
\[
\eta_m(x;P) 
\ = \  
\frac{f_m(x_i) - f_m(x_{i-1})}{\mu_C([x_{i-1},x_i])},
\]
where $x_{i-1} < x < x_i$ if $\mu_C([x_{i-1},x_i]) \neq 0$ and let
\[
\eta_m(x;P) = 0,
\]
where $x_{i-1} < x < x_i$ if $\mu_C([x_{i-1},x_i]) = 0$.
\end{x}
\vspace{0.25cm}

\begin{x}{\small\bf NOTATION} \ 
\[
\un{\eta}(x;P) 
\ = \  
(\eta_1(x;P), \ldots, \eta_M(x;P)).
\]
\\[-1cm]
\end{x}

\begin{x}{\small\bf THEOREM} \ 
\[
\int\limits_a^b \ 
\norm{\un{\theta}(x) - \un{\eta}(x;P)}^2 
\ \td\mu_C 
\ \leq \  
2
\hsx
\bigg[
\ell(C) \hsx - \hsx
\sum\limits_{i=1}^n \ 
\norm{\un{f}(x_i) - \un{f}(x_{i-1})}
\bigg].
\]
\\[-.75cm]

\un{Proof:}  
Given $P \in \sP[a,b]$, let $\sum^\prime$ denote a sum over intervals $[x_{i-1},x_i]$, where 
\\[-.5cm]

\noindent
$\norm{\un{\eta}(x;P)}^2 \neq 0$ and let $\sum^{\prime\prime}$ denote a sum over what remains.  
Now compute:

\allowdisplaybreaks
\begin{align*}
\int\limits_a^b \ 
\big|\hspace{-.05cm}\big|
\un{\theta}(x) &- \un{\eta}(x;P)
\big|\hspace{-.05cm}\big|^2 
\ \td\mu_C \ 
\\[18pt]
&= \ 
{\sum}^\prime \ 
\int\limits_{x_{i-1}}^{x_i} \ 
\norm{\un{\theta}(x) - \un{\eta}(x;P)}^2 
\ \td\mu_C 
\ + \  
{\sum}^{\prime\prime} \ 
\int\limits_{x_{i-1}}^{x_i} \ 
\norm{\un{\theta}(x)}^2 
\ \td\mu_C  
\\[18pt]
&= \ 
{\sum}^\prime \ 
\int\limits_{x_{i-1}}^{x_i} \ 
[\norm{\un{\theta}(x)}^2+\norm{\un{\eta}(x;P)}^2 - 2 \hsy \un{\theta}(x) \cdot \un{\eta}(x;P)] 
\ \td\mu_C
\\[18pt]
&\hspace{3cm} 
\ + \  
{\sum}^{\prime\prime} \int\limits_{x_{i-1}}^{x_i} \norm{\un{\theta}(x)}^2 
\ \td\mu_C
\\[18pt]
&= \ 
{\sum}^\prime \ 
\int\limits_{x_{i-1}}^{x_i} \ 
[1 +\norm{\un{\eta}(x;P)}^2 - 2\hsy \un{\theta}(x) \cdot \un{\eta}(x;P)] 
\ \td\mu_C
\ + \ 
{\sum}^{\prime\prime} \ 
\int\limits_{x_{i-1}}^{x_i} 1 
\ \td\mu_C
\\[18pt]
&= \ 
{\sum}^\prime \ 
[\mu_C([x_{i-1},x_i]) 
\ + \ 
\left[
\frac{\norm{\un{f}(x_i) - \un{f}(x_{i-1})}}{\mu_C([x_{i-1},x_i])}
\right]^2 
\ \mu_C([x_{i-1},x_i])
\\[18pt]
&\hspace{3cm} 
- 2 \hsy 
\frac{\norm{\un{f}(x_i) - \un{f}(x_{i-1})}^2}{\mu_C([x_{i-1},x_i])} 
\ + \ 
{\sum}^{\prime\prime} \ 
\mu_C([x_{i-1},x_i]) 
\\[18pt]
&\leq \ 
\ell(C) 
\ - \ 
{\sum}^\prime \ 
\frac{\norm{\un{f}(x_i) - \un{f}(x_{i-1})}^2}{\mu_C([x_{i-1},x_i])}
\\[18pt]
&\leq \ 
\ell(C) 
\ - \ 
{\sum}^\prime \ 
\norm{\un{f}(x_i) - \un{f}(x_{i-1})}
\\[18pt]
&\hspace{3cm}  
\ + \  
{\sum}^\prime \ 
\norm{\un{f}(x_i) - \un{f}(x_{i-1})} 
\left(1 - \frac{\norm{\un{f}(x_i) - \un{f}(x_{i-1})}}{\mu_C([x_{i-1},x_i])} \right)
\\[18pt]
&\leq \ 
\ell(C) 
\ - \ 
{\sum}^\prime \ 
\norm{\un{f}(x_i) - \un{f}(x_{i-1})}
\\[18pt]
&\qquad\qquad\qquad\qquad 
+ {\sum}^\prime \mu_C([x_{i-1},x_i])\left(1 - \frac{\norm{\un{f}(x_i) - \un{f}(x_{i-1})}}{\mu_C([x_{i-1},x_i])} \right)
\\[18pt]
&\leq \ 
\ell(C) - {\sum}^\prime \norm{\un{f}(x_i) - \un{f}(x_{i-1})}
\\[18pt]
&\hspace{3cm} 
\ + \  
{\sum}^\prime \left(\mu_C([x_{i-1},x_i]) - \norm{\un{f}(x_i) - \un{f}(x_{i-1})} \right)
\\[18pt]
&\leq  \ 
\ell(C) - {\sum}^\prime \norm{\un{f}(x_i) - \un{f}(x_{i-1})}
\\[18pt]
&\hspace{3cm} 
\ + \ 
{\sum}^\prime \ 
\mu_C([x_{i-1},x_i]) 
\ - \ 
{\sum}^\prime \ 
\norm{\un{f}(x_i) - \un{f}(x_{i-1})}
\\[18pt]
&=  \ 
\ell(C) 
\ + \ 
{\sum}^\prime \ 
\mu_C([x_{i-1},x_i]) - 2 \hsx {\sum}^\prime \  
\norm{\un{f}(x_i) - \un{f}(x_{i-1})}
\\[18pt]
&=   \ 
\ell(C)  
\ + \ 
\ell(C) - 2 \hsx {\sum}^\prime \ 
\norm{\un{f}(x_i) - \un{f}(x_{i-1})}
\\[18pt]
&=   \ 
2 \hsx \left[
\ell(C) 
\hsx - \hsx  
{\sum}^\prime \ 
\norm{\un{f}(x_i) - \un{f}(x_{i-1})}
\right]
\\[18pt]
&=   \ 
2 \hsx 
\bigg[
\ell(C)  
- \ 
\sum\limits_{i=1}^n \ 
\norm{\un{f}(x_i) - \un{f}(x_{i-1})}
\bigg].
\end{align*}
\\[-.75cm]
\end{x}

\begin{x}{\small\bf \un{N.B.}} \ 
By definition, $\mu_C([x_{i-1},x_i])$ is the length of the restriction of $C$ to $[x_{i-1},x_i]$, i.e.,
\[
\mu_C([x_{i-1},x_i]) 
\ = \  
s(x_i) - s(x_{i-1}).
\]
Moreover
\[
\norm{\un{f}(x_i) - \un{f}(x_{i-1})} 
\ \leq \  
s(x_i) - s(x_{i-1}).
\]
\\[-1cm]
\end{x}

So, if $\mu_C([x_{i-1},x_i]) = 0$, then
\[
\norm{\un{f}(x_i) - \un{f}(x_{i-1})} 
\ = \   
0
\quad \implies \quad \un{f}(x_i) = \un{f}(x_{i-1})
\]

\qquad
$\implies$
\begin{align*}
{\sum}^\prime \norm{\un{f}(x_i) - \un{f}(x_{i-1})} \ 
&= \ 
{\sum}^\prime \ 
\norm{\un{f}(x_i) - \un{f}(x_{i-1})}  
+ 
{\sum}^{\prime\prime} \ 
\norm{\un{f}(x_i) - \un{f}(x_{i-1})} 
\\[15pt]
&= \ 
\sum\limits_{i=1}^n \ 
\norm{\un{f}(x_i) - \un{f}(x_{i-1})}.
\end{align*}
Abbreviate
\[
\Lp^2 ([a,b], \tBO ([a,b]), \mu_C)
\]
to
\[
\Lp^2 (\mu_C).
\]
\\[-1cm]

\begin{x}{\small\bf APPLICATION} \ 
In $\Lp^2 (\mu_C)$,
\[
\lim\limits_{\norm{P} \rightarrow 0} \ 
\un{\eta} (-;P) = \un{\theta}.
\]
\\[-1cm]
\end{x}

\begin{x}{\small\bf SETUP} \ 
\\[-.5cm]

\qquad \textbullet \quad  
$C_0 \longleftrightarrow \un{f_0}:[a,b] \rightarrow \R^M$
\\[-.5cm]

\noindent
is a curve, continuous and rectifiable.
\\[-.25cm]

\qquad \textbullet \quad  
$C_k \longleftrightarrow \un{f_k}:[a,b] \rightarrow \R^M \qquad (k = 1, 2, \ldots)$
\\[-.5cm]

\noindent
is a sequence of curves, continuous and rectifiable.
\\[-.25cm]

Assumption:  $\un{f_k}$ converges uniformly to $\un{f_0}$ in $[a,b]$ and
\[
\lim\limits_{k \rightarrow \infty} \ 
\ell(C_k) 
\ = \  
\ell(C_0).
\]
\\[-1cm]
\end{x}

\begin{x}{\small\bf THEOREM} \ 
\[
\lim\limits_{\norm{Q} \rightarrow 0} \ 
\bigvee\limits_a^b \ (\un{f_k};Q) 
\ = \  
\ell(C_k) \qquad (Q \in \sP[a,b])
\]
uniformly in $k$, i.e., $\forall \ \varepsilon> 0, \ \exists \ \delta > 0$ such that
\[
\norm{Q} 
\ < \   
\delta \implies 
\bigg| \bigvee\limits_a^b \  (\un{f_k};Q) - \ell(C_k)
\bigg| 
\ < \   
\varepsilon
\]
for all $k = 1, 2, \ldots$, or still,
\[
\norm{Q} 
\ < \   
\delta \implies \ell(C_k) - \bigvee\limits_a^b \ (\un{f_k};Q) 
\ < \   
\varepsilon
\]
for all $k = 1, 2, \ldots \hsx .$ 
\\[-.25cm]
\end{x}

The proof will emerge in the lines to follow.  
Start the process by choosing $\delta_0 > 0$ such that
\[
\ell(C_0)  - \bigvee\limits_a^b \ (\un{f_0};P_0) 
\ < \   
\frac{\varepsilon}{4}
\]
provided $\norm{P_0} < \delta_0$.  Consider a $P \in \sP[a,b]:$
\[
a = x_0 < x_1 < \cdots < x_n = b
\]
with $\norm{P} < \delta_0$. 
Choose $\rho > 0$ such that
\[
\norm{\un{f_k}(c) - \un{f_k}(d)} 
\ < \   
\frac{\varepsilon}{4n} \qquad ([c,d] \subset [a,b])
\]
for all $k = 0, 1, 2, \ldots$, so long as $\abs{c - d} < \rho$ (equicontinuity).  
Take a partition $Q \in \sP[a,b]:$
\[
a = y_0 < y_1 < \cdots < y_m = b
\]
subject to
\[
\norm{Q} 
\ < \   
\gamma \equiv \min\limits_{i=1, \ldots,n} \ \bigg\{\rho, \frac{x_i - x_{i-1}}{2}\bigg\} 
\qquad (\implies \norm{Q} 
\ < \   
\delta_0).
\]
Put
\[
\sigma_k 
\ = \  
\sup\limits_{a \leq x \leq b} \ 
\norm{\un{f_k}(x) - \un{f_0}(x)}
\]
and let $k_0$ be such that
\[
k > k_0 \quad \implies \quad \sigma_k < \frac{\varepsilon}{4n} \quad \text{and} \quad \abs{\ell(C_k) - \ell(C_0)} < \frac{\varepsilon}{4}.
\]

The preparations complete, to minimize technicalities we shall suppose that each $I_j = [y_{j-1},y_j]$ is contained in just one $I_i = [x_{i-1},x_i]$ and write ${\sum}^{(i)}$ for a sum over all such $I_j$ $-$then
\begin{align*}
\bigvee\limits_a^b \ (\un{f_k};Q) \ 
&= \ 
\sum\limits_{j=1}^m \ 
v(\un{f_k};I_j) 
\\[15pt]
&= \ 
\sum\limits_{j=1}^m \ 
\norm{\un{f_k}(y_j) - \un{f_k}(y_{j-1})} 
\\[15pt]
&= \ 
\sum\limits_{i=1}^n \ 
{\sum}^{(i)} \norm{\un{f_k}(y_j) - \un{f_k}(y_{j-1})}
\\[15pt]
&\geq \ 
\sum\limits_{i=1}^n \ 
\norm{\un{f_k}(x_i) - \un{f_k}(x_{i-1})}.
\end{align*}
\vspace{0.2cm}

\begin{x}{\small\bf SUBLEMMA} \ 
Let $\un{A},\  \un{B},\  \un{C},\  \un{D} \in \R^M$ $-$then
\[
\norm{\un{C} - \un{D}} 
> 
\norm{\un{A} - \un{B}} - \norm{\un{A} - \un{C}} - \norm{\un{B} - \un{D}}.
\]
[In fact,
\begin{align*}
\norm{\un{A} - \un{B}} \ 
&= \ 
\norm{\un{A} - \un{C} + \un{C} - \un{D} + \un{D} - \un{B}}
\\[11pt] 
&\leq  \ 
\norm{\un{A} - \un{C}} + \norm{\un{C} - \un{D}} + \norm{\un{B} - \un{D}}\hsx ].
\end{align*}
\end{x}
\vspace{0.25cm}
Take
\[
\begin{cases}
\ \un{C} = \un{f_k}(x_i)
\\[4pt]
\ \un{D} = \un{f_k}(x_{i-1})
\end{cases}
\qquad
\begin{cases}
\ \un{A} = \un{f_0}(x_i)
\\[4pt]
\ \un{B} = \un{f_0}(x_{i-1})
\end{cases}
.
\]
Then
\begin{align*} 
 & \norm{\un{f_k}(x_i) - \un{f_k}(x_{i-1})} 
 \\[11pt]
&\hspace{1cm}
\geq  \
\norm{\un{f_0}(x_i) 
\ - \  
\un{f_0}(x_{i-1})} 
\ - \   
\norm{\un{f_0}(x_i) - \un{f_k}(x_{i})} 
\ - \  
\norm{\un{f_0}(x_{i-1}) - \un{f_k}(x_{i-1})},
\end{align*}
thus
\begin{align*} 
\sum\limits_{i=1}^n \ 
\norm{\un{f_k}(x_i) - \un{f_k}(x_{i-1})} \ 
&\geq \ 
\ell(C_0) - \frac{\varepsilon}{4} - n\sigma_k - n\sigma_k
\\[15pt]
&\geq \ 
\ell(C_0) - \frac{\varepsilon}{4} - \frac{\varepsilon}{4} - \frac{\varepsilon}{4}
\\[15pt]
&= \ 
\ell(C_0) - \frac{3 \varepsilon}{4}.
\end{align*}
But
\begin{align*}
k > k_0 
&\implies 
\abs{\ell(C_k - \ell(C_0} \ < \  \frac{\varepsilon}{4}
\\[11pt]
&\implies 
\ell(C_k) - \frac{\varepsilon}{4} \ < \   \ell(C_0).
\end{align*}
Therefore
\begin{align*}
\ell(C_0) - \frac{3 \varepsilon}{4} \ 
&> \ 
\ell(C_k) \hsx - \hsx  \frac{\varepsilon}{4} \hsx - \hsx \frac{3 \varepsilon}{4}
\\[11pt]
&= \ 
\ell(C_k) \hsx - \hsx \varepsilon.
\end{align*}
Thus:  $\forall \ k > k_0$,
\[
\ell(C_k) - \bigvee\limits_a^b \ (\un{f_k};Q) 
\ < \   
\varepsilon \qquad (\norm{Q} < \gamma).
\]
Finally, for $k \leq k_0$, let $\gamma_k$ be chosen so as to ensure that
\[
\ell(C_k) \hsx - \hsx \bigvee\limits_a^b \ (\un{f_k};Q) 
\ < \    
\varepsilon
\]
for all partitions $Q$ with $\norm{Q} < \gamma_k$.  Put now
\[
\delta 
= 
\min\limits_{1, \ldots, k_0} \{\gamma_1, \ldots \gamma_{k_0},\gamma\}.
\]
Then
\[
\norm{Q} 
\ < \   
\delta \quad \implies \quad \ell(C_k) - \bigvee\limits_a^b \ (\un{f_k};Q) 
\ < \    
\varepsilon
\]
for all $k = 1, 2, \ldots$ .
\\

Changing the notation (replace $Q$ by $P$), $\forall \ \varepsilon > 0$, $\exists \ \delta > 0$ such that
\[
\norm{P} 
\ < \   
\delta \quad \implies \quad \ell(C_k) - \bigvee\limits_a^b \ (\un{f_k};P) 
\ < \   
\varepsilon
\]
for all $k = 1, 2, \ldots$ .  
Consequently
\begin{align*}
\int\limits_a^b \ 
\norm{\un{\theta_k}(x) - \un{\eta_k}(x;P)}^2 
\ \td\mu_{C_k} \ 
&\leq \ 
2 \hsx
\bigg[
\ell(C_k) 
\hsx -  \hsx 
\sum\limits_{i=1}^n \ 
\norm{\un{f_k}(x_i) - \un{f_k}(x_{i-1})}
\bigg]
\\[15pt]
&= \ 
2 \hsx
\bigg[
\ell(C_k) 
\hsx -  \hsx  
\bigvee\limits_a^b \ (\un{f_k};P)
\bigg]
\\[15pt]
&< \ 
2 \hsx \varepsilon.
\end{align*}

\chapter{
$\boldsymbol{\S}$\textbf{8}.\quad LINE INTEGRALS}
\setlength\parindent{2em}
\setcounter{theoremn}{0}
\renewcommand{\thepage}{\S8-\arabic{page}}

\qquad
Let
\[
C \longleftrightarrow \un{f}:[a,b] \ra \R^M
\]
be a curve, continuous and rectifiable.
\\[-.25cm]

Suppose that
\[
F:[C] \times \R^M \ra \R,
\]
say
\[
F(\un{x},\un{t}) \quad (\un{x} \in [C], \un{t} \in \R^M).
\]
\\[-1cm]

\begin{x}{\small\bf DEFINITION} \ 
$F$ is a \un{parametric integrand} if $F$ is continuous in $(\un{x},\un{t})$ and $\forall \ K \geq 0$,
\[
F(\un{x},K\un{t}) \ = \  KF(\un{x},\un{t}).
\]
\\[-1.5cm]
\end{x}

\begin{x}{\small\bf EXAMPLE} \ 
Let
\[
F(\un{x},\un{t}) \ = \  (t_1^2 + \ldots + t_M^2)^\frac{1}{2}.
\]
\\[-1.5cm]
\end{x}

\begin{x}{\small\bf EXAMPLE} \ 
$(M = 2)$ Let
\[
F(x_1, x_2, t_1, t_2) \ = \  x_1 t_2 - x_2 t_1.
\]
\\[-1.5cm]
\end{x}

\begin{x}{\small\bf \un{N.B.}} \ 
If $F$ is a parametric integrand, then $\forall \ \un{x}$,
\[
F(\un{x},0) \ = \  0.
\]
\\[-1cm]
\end{x}

\begin{x}{\small\bf RAPPEL} \ 
\[
\norm{\un{\theta}} = 1
\]
almost everywhere.
\\[-.25cm]
\end{x}

\begin{x}{\small\bf LEMMA} \ 
Suppose that $F$ is a parametric integrand $-$then the integral
\[
I(C) 
\ \equiv \  
\int\limits_a^b F(\un{f}(x),\hsy \un{\theta}(x)) \ \td \mu_C
\]
exists.
\\[-.5cm]

PROOF \ \quad $[C] \times U(M)$ is a compact set on which $F$ is bounded.  
Since
\[
(\un{f}(x),\hsy \un{\theta}(x)) \in [C] \times U(M)
\]
almost everywhere, the function
\[
F(\un{f}(x),\hsy \un{\theta}(x)) 
\]
is Borel measurable and essentially bounded w.r.t the measure $\mu_C$.  Therefore
\[
I(C) \equiv \int\limits_a^b F(\un{f}(x),\hsy \un{\theta}(x)) \ \td \mu_C
\]
exists.
\\[-.25cm]

[Note: \  
The requirement ``homogeneous of degree 1'' in $t$ plays no role in the course of establishing the existence of $I(C)$.  
It will, however, be decisive in the considerations to follow.]
\\[-.25cm]
\end{x}

\quad Let $P \in \sP[a,b]$ and let $\xi_i$ be a point in $[x_{i-1},x_i]$ $(i = 1, \ldots, n)$.
\\[-.25cm]

\begin{x}{\small\bf THEOREM} \ 
If $F$ is a parametric integrand, then
\[
\lim\limits_{\norm{P} \ra 0} \ \sum\limits_{i=1}^n \ F(\un{f}(\xi_i), \un{f}(x_i)-\un{f}(x_{i-1}))
\]
exists and equals $I(C)$, denote it by the symbol
\[
\int\limits_C F, 
\]
and call it the \un{line integral} of $F$ along $C$.
\\[-.25cm]

PROOF \  
Fix $\varepsilon > 0$ and let $\tB(M)$ be the unit ball in $\R^M$.  Put
\[
M_F = \sup\limits_{[C] \times \tB(M)} \abs{F}.
\]
Choose $\gamma > 0:$
\[
\begin{cases}
\ \norm{\un{x_1} - \un{x_2}} < \gamma \qquad &(\un{x_1}, \un{x_2} \in [C])
\\[8pt]
\ \norm{\un{t_1} - \un{t_2}} < \gamma \qquad &(\un{t_1}, \un{t_2} \in \tB(M))
\end{cases}
\]

\qquad $\implies$
\[
\abs{F(\un{x_1},\un{t_1}) - F(\un{x_2},\un{t_2})} 
\ < \  
\frac{\varepsilon}{3\ell(C)}.
\]
Introduce $\un{\eta}(x;P)$ and set
\[
g(x;P) = F(\un{f}(\xi_i),\un{\eta}(x;P))
\]
if $x_{i-1} < x < x_i$ $-$then
\begin{align*}
\int\limits_a^b g(x;P) \ \td\mu_C \ 
&= \ 
\sum\limits_{i=1}^n \ 
F\left( \un{f}(\xi_i), \frac{\un{f}(x_i) - \un{f}(x_{i-1})}{\mu_C([x_{i-1},x_i])}\right) \mu_C([x_{i-1},x_i])
\\[15pt]
&= \ 
\sum\limits_{i=1}^n \ 
F( \un{f}(\xi_i),\un{f}(x_i) - \un{f}(x_{i-1}))
\end{align*}
modulo the usual convention if $\mu_C([x_{i-1},x_i]) = 0$.  Recall now that in $\Lp^2(\mu_C)$,
\[
\lim\limits_{\norm{P} \ra 0} \un{\eta}(-;P) = \un{\theta},
\]
hence $\un{\eta}(-;P)$ converges in measure to $\un{\theta}$, so there is a $\rho > 0$ such that for all $P$ with $\norm{P} < \rho$, 
\[
\norm{\un{\theta}(x) - \un{\eta}(x;P)} < \gamma
\]
except on a set $S_P$ of measure
\[
\mu_C(S_P) 
\ < \  
\frac{\varepsilon}{3M_F}.
\]
Define $\sigma:$
\[
\abs{t_1 - t_2} < \sigma \quad \implies \quad \norm{f(t_1) - f(t_2)} < \gamma.
\]
Let $\delta = \min(\sigma,\rho)$ and let $P$ be any partition with $\norm{P} < \delta$ $-$then
\begin{align*}
I(C) - \sum\limits_{i=1}^n \ F(\un{f}(\xi_i),&\un{f}(x_i) - \un{f}(x_{i-1}))  \ 
\\[15pt]
&= \ 
\int\limits_a^b  \ 
F(\un{f}(x),\un{\theta}(x)) \ \td\mu_C  
- 
\int\limits_a^b \ g(x;P)  \ \td\mu_C
\\[15pt]
&= \ 
\int\limits_a^b  \ 
[F(\un{f}(x),\un{\theta}(x)) -g(x;P)] \ \td\mu_C.
\end{align*}
By definition, $\delta \leq \rho$, hence
\[
\norm{\un{\theta}(x) - \un{\eta}(x;P)} 
\ < \  
\gamma
\]
except in $S_P$, and
\[
\norm{\un{f}(x) - \un{f}(\xi_i)} \ < \  \gamma
\]
since
\[
\abs{x - \xi_i} < \gamma \quad (x_{i-1} \leq x \leq x_i).
\]
To complete the argument, take absolute values:

\begin{align*}
\bigg|
I(C) \hsx - \hsx \sum\limits_{i=1}^n \ 
F(\un{f}(\xi_i),&\hsx \un{f}(x_i) - \un{f}(x_{i-1}))
\bigg|
\\[15pt]
&\leq \ 
\int\limits_a^b  \ 
\abs{F(\un{f}(x),\un{\theta}(x)) - g(x;P)} 
\ \td\mu_C
\\[15pt]
&= \ 
\int\limits_{[a,b] - S_P}  \ 
\abs{\ldots} 
\ \td\mu_C
+ 
\int\limits_{S_P} \ 
\abs{\ldots} 
\ \td\mu_C.
\end{align*}
\\[-1cm]

\qquad \textbullet \quad On $[a,b]$ $-$ $S_P$ at an index $i$, 

\begin{align*}
\abs{F(\un{f}(x),\un{\theta}(x)) - g(x;P)} \ 
&= \ 
\abs{F(\un{f}(x),\hsy \un{\theta}(x))  - F(\un{f}(\xi_i),\hsy \un{\eta}(x;P))}
\\[15pt]
&\leq \ 
\frac{\varepsilon}{3 \hsy \ell(C)}.
\end{align*}
Here, of course, up to a set of measure 0, 
\[
\un{\theta}(x) \in \tB(M) \text{ and } \eta (x;P) \in \tB(M).
\]
Therefore
\begin{align*}
\int\limits_{[a,b] - S_P}  \ 
\abs{\ldots} 
\ \td\mu_C \  
&\leq \ 
\frac{\varepsilon}{3 \hsy \ell(C)}
\\[15pt]
&= \  \frac{\varepsilon}{3}.
\end{align*}

\qquad \textbullet \quad On $S_P$,

\[
\begin{cases}
\ \abs{F(\un{f}(x),\un{\theta}(x)) } \hsx \leq  \hsx M_F
\\[8pt]
\ \abs{F(\un{f}(\xi_i),\un{\eta}(x;P))} \hsx \leq  \hsx M_F
\end{cases}
.
\\
\]
Therefore
\begin{align*}
\int\limits_{S_P}  \ 
\abs{\ldots} 
\ \td\mu_C \  
&\leq \ 
2M_F \int\limits_{S_P} \ 1 
\ \td \mu_C
\\[15pt]
&= \ 
2 \hsx M_F \hsx \mu_C(S_P) 
\\[15pt]
&< \ 
2 \hsx M_F \hsx \frac{\varepsilon}{3M_F}
\\[15pt]
&= \ 
\frac{2\varepsilon}{3}.
\end{align*}
So in conclusion,

\begin{align*}
\int\limits_{[a,b] - S_P}  \ 
\abs{\ldots}
\ \td\mu_C 
\hsx + \hsx 
\int\limits_{S_P} \ 
\abs{\ldots}
\ \td\mu_C \ 
&< \ \frac{\varepsilon}{3} + \frac{2\varepsilon}{3}
\\[15pt]
&= \
\varepsilon \qquad (\norm{P} < \delta)
\end{align*}
and
\[
I(C) \ = \  \int\limits_C F.
\]
\\[-1cm]
\end{x}

\begin{x}{\small\bf \un{N.B.}} \ 
The end result is independent of the choice of the $\xi_i$.
\\[-.25cm]
\end{x}

\begin{x}{\small\bf THEOREM} \ 
If $f_1, \ldots, f_M \in \AC[a,b]$, then for any parametric \text{integrand $F$,}
\[
\int\limits_C F  \ = \  \int\limits_a^b F(f_1(x), \ldots, f_M(x), f_1^\prime(x), \ldots, f_M^\prime(x)) \ \td x,
\]
the integral on the right being in the sense of Lebesgue.
\\[-.25cm]

PROOF \  
The absolute continuity of the $f_m$ implies that
\[
\mu_C([c,d]) \ = \  \int\limits_c^d \norm{f^\prime}\ \td x
\]
for every subinterval $[c,d] \subset [a,b]$, hence $\mu_C$ is absolutely continuous w.r.t Lebesgue measure.  
It is also true that $\nu_m$ is absolutely continuous w.r.t Lebesgue measure.  
This said, write
\[
f_m^\prime 
\ = \  
\frac{\td f_m}{\td x} 
\ = \  
\frac{\td\nu_m}{\td x} 
\ = \  
\frac{\td\nu_m}{\td\mu_C} \frac{\td\mu_C}{\td x} 
\ = \  
\theta_m\frac{\td\mu_C}{\td x}.
\]
Then
\begin{align*}
I(C) 
&= \ 
\int\limits_a^b \ 
F(\un{f}(x),\hsy \un{\theta}(x)) \ \td\mu_C
\\[15pt]
&= \ 
\int\limits_a^b \ 
F(\un{f}(x),\hsy \un{\theta}(x)) \frac{\td\mu_C}{\td x} \td x
\\[15pt]
&= \ 
\int\limits_a^b \ 
F\left(\un{f}(x),\hsy \un{\theta}(x) \frac{\td\mu_C}{\td x} \right)\ \td x,
\end{align*}
where
\[
\frac{\td\mu_C}{\td x} 
\ = \  
\norm{\un{f}^\prime} 
\ \geq \ 
 0.
\]
Continuing
\begin{align*}
I(C) \ 
&= \ 
\int\limits_a^b \ 
F\left(f_1(x), \ldots, f_M(x),\hsy \theta_1(x)\frac{\td\mu_C}{\td x}, \ldots, \theta_M(x)\frac{\td\mu_C}{\td x}\right) \td x
\\[15pt]
&= \ 
\int\limits_a^b \ 
F(f_1(x), \ldots, f_M(x), f_1^\prime(x), \ldots, f_M^\prime(x)) \ \td x, 
\end{align*}
the integrals being in the sense of Lebesgue.
\\[-.25cm]
\end{x}

Let
\[
\begin{cases}
\ C \longleftrightarrow \un{f}:[a,b] \ra \R
\\
\ D \longleftrightarrow \un{g}:[a,b] \ra \R
\end{cases}
\]
\\[-.75cm]

\noindent
be curves, continuous and rectifiable.
\\[-.25cm]

\begin{x}{\small\bf RAPPEL} \ 
If $C$ and $D$ are Fr\'echet equivalent, then
\[
[C] \ = \  [D] 
\quad \text{and} \quad 
\ell(C) \ = \  \ell(D).
\]
\\[-1.25cm]
\end{x}

\begin{x}{\small\bf THEOREM} \ 
If $C$ and $D$ are Fr\'echet equivalent and if $F$ is a parametric integrand, then
\[
\int\limits_C F = \int\limits_D F.
\]
\\[-1cm]

PROOF \ \ 
Fix $\varepsilon >0$ and choose $\delta > 0:$
\\[-.25cm]

\qquad \textbullet \quad $P \in \sP[a,b]$ \ \& \  $\norm{P} < \delta$ $\implies$

\[
\bigg|
I(C) - \sum\limits_{i=1}^n F(\un{f}(\xi_i),\un{f}(x_i) - \un{f}(x_{i-1}))
\bigg|
\ < \ 
\frac{\varepsilon}{3}.
\]

\qquad \textbullet \quad $Q \in \sP[c,d]$ \ \& \  $\norm{Q} < \delta$ $\implies$

\[
\bigg|
I(D) - \sum\limits_{j=1}^m F(\un{f}(\xi_j),\un{f}(y_j) - \un{f}(y_{j-1}))
\bigg| 
\ < \ 
\frac{\varepsilon}{3}.
\]

\noindent
Fix $P$ and $Q$ satisfying these conditions and let $k$ be the number of intervals in $P$ 
and let $\ell$ be the number of intervals in $Q$.  
Fix $\gamma > 0$ such that
\[
\abs{F(\un{x}_1,\un{t}_1) - F(\un{x}_2,\un{t}_2)} 
\ < \  
\frac{\varepsilon}{3(k+\ell)}
\]
when
\[
\norm{\un{x}_1 - \un{x}_2} 
\ < \  
\gamma \qquad (\un{x}_1, \un{x}_2 \in [C] = [D])
\]
and
\[
\norm{\un{t}_1 - \un{t}_2} 
\ < \  
2 \gamma \qquad (\norm{\un{t}_1} 
\ \leq \  
\ell(C), \  \norm{\un{t}_2} \leq \ell(D)).
\]
Let $\phi:[a,b] \ra [c,d]$ be a homeomorphism $(\phi(a) = c, \ \phi(b) = d)$ such that
\[
\norm{\un{f}(x) - \un{g}(\phi(x))} 
\ < \  
\gamma \qquad (x \in [a,b]).
\]
Let
\[
P^*: \ a = x_0^* < x_1^* < \cdots < x_r^* = b
\]
be the partition obtained from $P$ by adjoining the images under $\phi^{-1}$ of the partition points of $Q$.  
Let
\[
Q^*: \ c = y_0^* < y_1^* < \cdots < y_s^* = d
\]
be the partition obtained from $Q$ by adjoining the images under $\phi$ of the partition points of $P$.  
So, by construction, $r = s$, either one is $\leq k + \ell$, and $y_p^* = \phi(x_p^*) \ (p = 0, 1, \ldots, q)$.  
Choose a point $\xi_p \in [x_{p-1}^*, x_p^*]$ and work with 
\[
\un{f}(\xi_p ) \quad \text{and} \quad \un{g}(\phi(\xi_p )).
\]
Then
\begin{align*}
\abs{I(C) - I(D)} \ 
&\leq \ 
\bigg|
I(C) - \sum\limits_{p=1}^q \ 
F(\un{f}(\xi_p),\un{f}(x_p^*)  - \un{f}(x_{p-1}^*))
\bigg|
\\[15pt]
&\qquad + 
\sum\limits_{p=1}^q \ 
\abs{F(\un{f}(\xi_p),\un{f}(x_p^*)  - \un{f}(x_{p-1}^*)) - F(\un{g}(\phi(\xi_p)),\un{g}(y_p^*)  - \un{g}(y_{p-1}^*)) }
\\[15pt]
&\qquad \qquad + 
\bigg|
\sum\limits_{p=1}^q \ 
F(\un{g}(\phi(\xi_p)),\un{g}(y_p^*)  - \un{g}(y_{p-1}^*)) - I(D)
\bigg|.
\end{align*}
Since
\[
\begin{cases}
\ \norm{P^*} \ \leq \ \norm{P} \ < \ \delta
\\[8pt]
\ \norm{Q^*} \ \leq \ \norm{Q} \ < \ \delta
\end{cases}
,
\]
the first and third terms are each $\ds < \frac{\varepsilon}{3}$.  
As for the middle term,
\[
\norm{\un{f}(\xi_p) - \un{g}(\phi(\xi_p))} < \gamma
\]
and
\allowdisplaybreaks
\begin{align*}
&\norm{\un{f}(x_p^*) - \un{f}(x_{p-1}^*)  - \un{g}(y_p^*) + \un{g}(y_{p-1}^*)} \ 
\\[11pt]
&\hspace{2cm}\leq \
\norm{\un{f}(x_p^*) - \un{g}(y_p^*)} + \norm{\un{f}(x_{p-1}^*) - \un{g}(y_{p-1}^*)}
\\[11pt]
&\hspace{2cm}= \ 
\norm{\un{f}(x_p^*) - \un{g}(\phi(x_p^*))} + \norm{\un{f}(x_{p-1}^*) - \un{g}(\phi(x_{p-1}^*))}
\\[11pt]
&\hspace{2cm}< \ 
\gamma + \gamma
\\[11pt]
&\hspace{2cm}= \ 
2\gamma.
\end{align*}
Therefore the middle term is
\[
< \  q \hsx \frac{\varepsilon}{3(k+\ell)} 
\ = \ 
\frac{q}{k+\ell} \ \frac{\varepsilon}{3} 
\ < \ 
\frac{\varepsilon}{3}.
\]
And finally
\[
\abs{I(C) - I(D)} 
\ < \  
\frac{\varepsilon}{3}+ \frac{\varepsilon}{3} + \frac{\varepsilon}{3} 
\ = \ 
\varepsilon
\]

\qquad $\implies$
\[
I(C) \ = \  I(D) \qquad (\varepsilon \downarrow 0)
\]

\qquad $\implies$
\[
\int\limits_C \ F \ = \  \int\limits_D \ F.
\]
\\[-1cm]
\end{x}

\begin{x}{\small\bf SETUP} \ 
\\[-.25cm]

\qquad \textbullet \quad 
$
C_0 \longleftrightarrow \un{f_0} :[a,b] \ra \R^M
$
\\[-.5cm]

\noindent
is a curve, continuous and rectifiable.
\\[-.25cm]

\qquad \textbullet \quad 
$
C_k \longleftrightarrow \un{f_k}:[a,b] \ra \R^M \qquad (k = 1, 2, \ldots)
$
\\[-.5cm]

\noindent
is a sequence of curves, continuous and rectifiable.
\\[-.25cm]

\quad Assumption:  $\un{f_k}$ converges uniformly to $\un{f_0}$ in $[a,b]$ and
\[
\lim\limits_{k \ra \infty} \ell(C_k) 
= \   
\ell(C_0).
\]
\\[-1cm]
\end{x}

\begin{x}{\small\bf THEOREM} \ 

\[
\lim\limits_{k \ra \infty} \ I(C_k) 
\ = \  
I(C_0) 
\]
or still,

\[
\lim\limits_{k \ra \infty} \ \int\limits_{C_k} \ F 
\ = \ 
\int\limits_{C_0} \ F.
\]
\end{x}

\chapter{
$\boldsymbol{\S}$\textbf{9}.\quad QUASI ADDITIVITY}
\setlength\parindent{2em}
\setcounter{theoremn}{0}
\renewcommand{\thepage}{\S9-\arabic{page}}

\begin{x}{\small\bf DATA} \ 
$A$ is a nonempty set, $\sI = \{I\}$ is a nonempty collection of subsets of $A$, 
$\sD = \{D\}$ is a nonempty collection of nonempty finite collections of 
$D = [I]$ of sets $I \in \sI$, and $\delta$ is a real valued function defined on $\sD$.
\\[-.25cm]
\end{x}

\begin{x}{\small\bf DEFINITIONS} \ 
The sets $I \in \sI$ are called \un{intervals}, 
the collections $D \in \sD$ are called \un{systems}, and the function $\delta$ is called a \un{mesh}.
\\[-.25cm]
\end{x}

\begin{x}{\small\bf ASSUMPTIONS} \ 
$A$ is a nonempty topological space, each interval $I$ has a nonempty interior, 
the intervals of each system $D$ are nonoverlapping$:$ $I_1, I_2 \in D$, $I_1 \neq I_2$ 
\[
\implies \qquad
\begin{cases}
\ \text{int} \ I_1 \cap \  c\ell I_2 \ = \ \emptyset 
\\[4pt]
\ c\ell  \ I_1 \cap \text{int}\  I_2 \ = \ \emptyset 
\end{cases}
.
\]
\\[-.75cm]
\end{x}

\begin{x}{\small\bf ASSUMPTION} \ 
For each system $D$, $0 < \delta(D) < +\infty$, and each $\varepsilon > 0$, 
there are systems with $\delta(D) < \varepsilon$.
\\[-.25cm]
\end{x}

\begin{x}{\small\bf REMARK} \ 
In the presence of $\delta$, one is able to convert $\sD$ into a directed set with direction ``$\gg$'' 
by defining $D_2 \gg D_1$ iff $\delta(D_2) < \delta(D_1)$.
\\[-.25cm]
\end{x}

\begin{x}{\small\bf EXAMPLE} \ 
Take $A = [a,b]$ and let $\sI = \{I\}$ be the collection of all closed subintervals of $[a,b]$.  
Take for $\sD$ the class of all partitions $D$ of $[a,b]$, 
i.e., $\sD = \sP[a,b]$, and let $\delta(D)$ be the norm of $D$.
\\[-.5cm]

[Note: \    
Strictly speaking, an element of $\sP[a,b]$ is a finite set $P = \{x_0, \ldots, x_n\}$, 
where
\[
a = x_0 < x_1 < \cdots < x_n = b,
\]
the associated element $D$ in $\sD$ being the set
\[
[x_{i-1},x_i] \qquad (i = 1, \ldots, n).]
\]
\\[-1cm]
\end{x}

\begin{x}{\small\bf DEFINITION} \ 
An \un{interval function} is a function $\phi: \sI \ra \R^M$.
\\[-.5cm]

[Note: \   Associated with $\phi$ is the interval function $\norm{\phi}$, as well as the 
\[
\phi_m, \ \abs{\phi_m}, \ 
\begin{cases}
\ \phi_m^+
\\[4pt]
\ \phi_m^-
\end{cases}
\qquad (m = 1, \ldots, M).]
\]
\\[-1cm]
\end{x}

\begin{x}{\small\bf NOTATION} \ 
Given an interval function $\phi$, a subset $S \subset A$, and a system $D = [I]$, put

\[
\sum \ [\phi,S,D] 
\ = \  
\sum\limits_I s(I,S) \  \phi(I),
\]
where $\ds\sum\limits_I$ ranges over all $I \in D$ and $s(I,S) = 1$ or 0 depending on whether $I \subset S$ or $I \not\subset S$.
\\[-.5cm]

[Note: \   Take for $S$ the empty set $\emptyset$ $-$then $I \subset \emptyset$ is inadmissible 
($I$ has nonempty interior) and $I \not\subset S$ gives rise to zero.  Therefore
\[
\sum \ [\phi,\emptyset,D] 
\ = \  
0.]
\]
\\[-1cm]
\end{x}

\begin{x}{\small\bf \un{N.B.}} \ 
The \un{absolute situation} is when $S = A$, thus  in this case, 
\[
\sum \ [\phi,A,D] 
\ \equiv \  
\sum \ [ \phi,D] \ = \  \sum\limits_I \ \phi(I).
\]
\\[-1cm]
\end{x}

\begin{x}{\small\bf DEFINITION} \ 
Given an interval function $\phi$ and a subset $S \subset A$, the \un{BC-integral} of $\phi$ over $S$ is
\[
\lim\limits_{\delta(D) \ra 0} \  \sum \ [\phi,S,D]
\]
provided the limit exists in $\R^M$.
\\[-.5cm]

[Note: \   
B = Burkill and C = Cesari.]
\\[-.25cm]
\end{x}

\begin{x}{\small\bf NOTATION} \ 
The BC-integral of $\phi$ over $S$ is denoted by
\[
\tBC \int\limits_S \phi.
\]
\\[-1cm]
\end{x}

\begin{x}{\small\bf EXAMPLE} \ 
\[
\tBC \int\limits_\emptyset \ \phi 
\ = \  
\un{0} \quad (\in \R^M).
\]
\\[-1cm]
\end{x}

\begin{x}{\small\bf DEFINITION} \ 
An interval function $\phi$ is \un{quasi additive} on $S$ if for each $\varepsilon > 0$ there exists 
$\eta(\varepsilon,S) > 0$ 
such that if 
$D_0 = [I_0]$ is any system subject to $\delta(D_0) < \eta(\varepsilon,S)$, 
there also exists 
$\lambda(\varepsilon,S,D_0) > 0$ such that for every system $D = [I]$ with 
$\delta(D) < \lambda(\varepsilon,S,D_0)$, 
the relations
\begin{align*}
&(\tqa_1 - S) \hsx 
\sum\limits_{I_0} s(I_0,S) \ 
\normx{\sum\limits_I s(I,I_0), \phi(I) - \phi(I_0)} \ < \ \varepsilon
\\[15pt]
&(\tqa_2 - S) \hsx 
\sum\limits_{I} s(I,S) \ 
\bigg[
1 \hsx - \hsx \sum\limits_{I_0} s(I,I_0), s(I_0,S)
\bigg] 
\norm{\phi(I)} 
\ < \ \varepsilon
\end{align*}
obtain.
\\[-.25cm]
\end{x}

\begin{x}{\small\bf \un{N.B.}} \ 
In the absolute situation, matters read as follows:  
An interval function $\phi$ is \un{quasi additive} if for each $\varepsilon > 0$  
there exists $\eta(\varepsilon) > 0$ such that if $D_0 = [I_0]$ is any system subject to 
$\delta(D_0) < \eta(\varepsilon)$ there  exists $\lambda(\varepsilon,D_0) > 0$ 
such that for every system $D = [I]$ with $\delta(D) < \lambda(\varepsilon,D_0)$, the relations
\begin{align*}
&(\tqa_1 - A) \ 
\sum\limits_{I_0} \ 
\normx{\sum\limits_{I \subset I_0} \ 
\phi(I) - \phi(I_0)} 
\ < \ 
\varepsilon 
\\[15pt]
&(\tqa_2 - A) \
\sum\limits_{I \not\subset I_0} \ 
\norm{\phi(I)} 
\ < \ 
\varepsilon
\end{align*}
obtain.
\\[-.5cm]

[Note: \   
The sum
\[
\sum\limits_{I \not\subset I_0} \norm{\phi(I)} 
\]
is over all $I \in D$, $I \not\subset I_0$ for any $I_0 \in D_0.]$
\\[-.25cm]
\end{x}

So, under the preceding conditions,
\[
\sum\limits_I \ \phi(I) - \sum\limits_{I_0} \ \phi(I_0) 
\ = \  
\sum\limits_{I_0} \ 
\bigg[
\sum\limits_{I \subset I_0} \ \phi(I) - \phi(I_0)
\bigg] 
+ 
\sum\limits_{I \not\subset I_0} \ \phi(I)
\]
\\[-1cm]

$\implies$
\[
\normxx{\sum\limits_I \ \phi(I) - \sum\limits_{I_0} \ \phi(I_0)} 
\ < \  2 
\varepsilon.
\]
\\[-.5cm]
 
\begin{x}{\small\bf THEOREM} \ 
If $\phi$ is quasi additive on $S$, then
\[
\tBC \int\limits_S \phi
\]
exists.
\\[-.25cm]

PROOF \ 
To simplify the combinatorics, take $S = A$.  
Given $\varepsilon > 0$, let 
$\eta(\varepsilon)$, $D_0$, 
$\lambda(\varepsilon,D_0)$ 
be per $\tqa_1$-$A$, $\tqa_2$-$A$ and suppose that $D_1, \ D_2 \in \sD$, where
\[
\begin{cases}
\ \delta(D_1) \ < \ \lambda(\varepsilon,D_0) 
\\[4pt]
\ \delta(D_2) \ < \ \lambda(\varepsilon,D_0)
\end{cases}
.
\]
 
\noindent
Then
\[
\begin{cases}
\ds\ \normx{\sum\limits_{I_1} \phi(I_1) -  \sum\limits_{I_0} \phi(I_0)} 
\ < \ 
2 \hsy \varepsilon  
\\[26pt]
\ds\ \normx{\sum\limits_{I_2} \phi(I_2) -  \sum\limits_{I_0} \phi(I_0)} 
\ < \ 
2 \hsy \varepsilon 
\end{cases}
\]
\\[-.5cm]

$\implies$
\[
\normx{\sum\limits_{I_1} \phi(I_1) -  \sum\limits_{I_2} \phi(I_2)} 
\ < \ 
4 \hsy \varepsilon.
\]
\\[-.5cm]

\noindent
Therefore $\ds\tBC \int\limits_A \phi$ exists.
\\
\end{x}

\begin{x}{\small\bf REMARK} \ 
\\[-.5cm]

\qquad \textbullet \quad If the $\phi_m$ $(m = 1, \ldots, M)$ are quasi additive, then $\phi$ is quasi additive.
\\

\qquad \textbullet \quad If the $\abs{\phi_m}$ $(m = 1, \ldots, M)$ are quasi additive, then $\norm{\phi}$ is quasi additive.
\\[-.25cm]
\end{x}

\begin{x}{\small\bf DEFINITION} \ 
$A$ real valued interval function $\psi$ is \un{quasi subadditive} on $S$ if for each $\varepsilon > 0$ 
there exists $\eta(\varepsilon,S) > 0$ such that if $D_0 = [I_0]$ is any system subject to 
$\delta(D_0) < \eta(\varepsilon,S)$  
there also exists $\lambda(\varepsilon,S,D_0) > 0$ such that 
for every system $D = [I]$ with $\delta(D) < \lambda(\varepsilon,S,D_0)$ the relation
\[
(\tqsa S) \ 
\sum\limits_{I_0} \ s(I_0,S) \ 
\bigg[
\sum\limits_I \ s(I,I_0) \psi(I) - \psi(I_0)
\bigg]^- 
\ < \ 
\varepsilon
\]
obtains.
\\[-.25cm]
\end{x}

\begin{x}{\small\bf \un{N.B.}} \ 
In the absolute situation, matters read as follows$:$ \ldots
\[
(\tqsa A) \ 
\sum\limits_{I_0} \ 
\bigg[\sum\limits_{I \subset I_0} \ \psi(I) - \psi(I_0)\bigg]^- 
\ < \ 
\varepsilon.
\]
\\[-1cm]
\end{x}

\begin{x}{\small\bf LEMMA} \ 
If $\psi: \sD \ra \R_{\geq 0}$ is nonnegative and quasi subadditive on  $S$, then
\[
\tBC \int\limits_S \psi
\]
exists $( +\infty$ is a permissible value$)$.
\\[-.25cm]
\end{x}

\begin{x}{\small\bf THEOREM} \ 
If $\psi: \sI \ra \R_{\geq 0}$ is nonnegative and quasi subadditive on $S$ and if 
\[
\tBC \int\limits_S \psi
\]
is finite, then $\psi$ is quasi additive on  $S$ .
\\[-.5cm]

PROOF \ 
To simplify the combinatorics, take $S = A$.  Since
\[
\tBC \ \int\limits_A \psi
\]
exists and is finite, given $\varepsilon > 0$ there is a number $\mu(\varepsilon) > 0$ such that for any $D_0 = [I_0] \in \sD$ with $\delta(D_0) < \mu(\varepsilon)$, we have
\[
\bigg| 
\tBC \ \int\limits_A \ \psi - \sum\limits_{I_0} \psi(I_0)
\bigg| 
\ < \ 
\frac{\varepsilon}{3},
\]
where $\ds\sum\limits_{I_0}$ is a sum ranging over all $I_0 \in D_0$.  Now choose $D_0$ in such a way that
\[
\delta(D_0) 
\ < \  
\min\{ \mu(\varepsilon), \eta(\varepsilon/6) \},
\]
take
\[
\lambda^\prime(\varepsilon)  
\ = \  
\min\{ \mu(\varepsilon), \lambda(\varepsilon/6, D_0) \},
\]
and consider any system $D = [I]$ with $\delta(D) < \lambda^\prime$.  Since $\psi$ is quasi subadditive, 
\[
\sum\limits_{I_0} 
\bigg[ 
\sum\limits_{I \subset I_0} \ 
\psi(I) - \psi(I_0)
\bigg]^- 
\ < \ 
\frac{\varepsilon}{6}.
\]
On the other hand,
\[
\bigg|
\tBC \int\limits_A \ \psi - \sum\limits_{I} \psi(I)
\bigg|
\ < \ 
\frac{\varepsilon}{3}.
\]
Denote by $\ds\sum^\prime$ a sum over all $I \in D$ with $I \not\subset I_0$ for any $I_0 \in D_0$ $-$ then
 
 \allowdisplaybreaks
\begin{align*}
0 \ 
&\leq \ 
\sum\limits_{I_0} \ 
\bigg|
\sum\limits_{I \subset I_0} \ 
\psi(I) - \psi(I_0)
\bigg| 
\ + \ \sum{}^\prime \ \psi(I)
\\[15pt]
&= \ 
\sum\limits_{I_0} \ 
\bigg[ 
\sum\limits_{I \subset I_0} \ 
\psi(I) - \psi(I_0)
\bigg] 
+ 2 \hsx  
\sum\limits_{I_0} 
\bigg[ 
\sum\limits_{I \subset I_0} \psi(I) - \psi(I_0)
\bigg]^- 
\\[15pt]
&\hspace{2cm}
+ \sum{}^\prime \psi(I)
\\[15pt]
&= \ 
\bigg[ 
\sum\limits_I \psi(I) - \tBC\int\limits_A \psi
\bigg] 
- 
\bigg[
\sum\limits_{I_0} \ 
\psi(I_0) - \tBC\int\limits_A \psi 
\bigg] 
\\[15pt]
&\hspace{2cm}
+ 
2 \hsx 
\sum\limits_{I_0} \hsx 
\bigg[ 
\sum\limits_{I \subset I_0} \psi(I) - \psi(I_0)
\bigg]^-
\\[15pt]
&\leq \ 
\frac{\varepsilon}{3} \ + \  \frac{\varepsilon}{3} \ + \  2 \hsx \frac{\varepsilon}{6}
\\[15pt]
&= \ 
\varepsilon.
\end{align*}
The requirements for quasi additivity are thus met.
\\[-.25cm]
\end{x}

\begin{x}{\small\bf THEOREM} \ 
Suppose that $\phi:\sI \ra \R^M$ is quasi additive on $S$ 
$-$then $\norm{\phi} :\sI \ra \R_{\geq 0}$ is quasi additive on $S$.
\\[-.25cm]

PROOF \ 
Fix $\varepsilon > 0$, take $S = A$, and in the notation above, introduce $\eta(\varepsilon)$, $D_0 = [I_0]$, $\lambda(\varepsilon,D_0)$, $D = [I]$ $-$then the objective is to show that
\[
\sum\limits_{I_0}  \ 
\bigg[ \sum\limits_{I \subset I_0} \ 
\norm{\phi(I)} - \norm{\phi(I_0)} 
\bigg]^- 
\ < \ 
\varepsilon.
\]
To this end, let
\[
\Phi(I_0) 
\ = \  
\sum\limits_{I \subset I_0} \ \phi(I) - \phi(I_0).
\]
Then
\allowdisplaybreaks
\begin{align*}
\norm{\phi(I_0) + \Phi(I_0)} 
&= \ 
\sum\limits_{I \subset I_0} \ 
\phi(I) \ 
\\[15pt]
&= \ 
\bigg[  
\sum\limits_{m=1}^M \ 
\bigg( 
\sum\limits_{I \subset I_0} \ 
\phi_m(I)
\bigg)^2\ 
\bigg]^\frac{1}{2}
\\[15pt]
&\leq\ 
\sum\limits_{I \subset I_0}  \ 
\bigg[  
\sum\limits_{m=1}^M \ 
\phi_m(I)
^2\ 
\bigg]^\frac{1}{2}
\\[15pt]
&=\ 
\sum\limits_{I \subset I_0} \ 
\norm{\phi(I)}.
\end{align*}
Meanwhile
\[
\phi(I_0) = [\phi(I_0) + \Phi(I_0)] + [-\Phi(I_0)]
\]

$\implies$
\begin{align*}
\sum\limits_{I \subset I_0} \ 
\norm{\phi(I)} - \norm{\phi(I_0)} \ 
&\geq \ 
\norm{\phi(I_0) + \Phi(I_0)} - \norm{\phi(I_0)} 
\\[15pt]
&\geq -\norm{\Phi(I_0)}
\end{align*}

$\implies$
\[
\bigg[ 
\sum\limits_{I \subset I_0} \norm{\phi(I)} - \norm{\phi(I_0)}
\bigg]^- 
\ \leq \ 
\norm{\Phi(I_0)} 
\]

$\implies$
\begin{align*}
\sum\limits_{I_0} \ 
\bigg[ 
\sum\limits_{I \subset I_0} \ 
\norm{\phi(I)} - \norm{\phi(I_0)}
\bigg]^- \ 
&\leq \ 
\sum\limits_{I_0} \ 
\norm{\Phi(I_0)} 
\\[15pt]
&= \ 
\sum\limits_{I_0} \ 
\normx{\sum\limits_{I \subset I_0} \ {\phi(I)} - {\phi(I_0)}}
\\[15pt]
&< \ 
\varepsilon,
\end{align*}
$\phi$ being quasi additive.
\\[-.25cm]
\end{x}

\begin{x}{\small\bf APPLICATION} \ 
If $\phi:\sI \ra \R^M$ is quasi additive, then the interval functions
\[
I \ra \abs{\phi_m(I)} \quad (m = 1, \ldots, M)
\]
are quasi subadditive.
\\[-.5cm]

[In fact, the quasi additivity of the $\phi$ implies the quasi additivity of the $\phi_m$ and
\[
\norm{\phi_m} = \abs{\phi_m}.]
\] 

[Note: \   It is also true that $\phi_m^+, \ \phi_m^-$ are quasi subadditive.]
\\[-.25cm]
\end{x}

\begin{x}{\small\bf LEMMA} \ 
If $\phi:\sI \ra \R^M$ is quasi additive on $S$ and if
\[
\tBC \int\limits_S \ \norm{\phi} 
\ < \ 
+\infty,
\] 
then $\phi$ is quasi additive on every subset $S^\prime \subset S$.
\\

PROOF   \
First of all, $\norm{\phi}$ is quasi subadditive on  $S$ , hence also on $S^\prime$.  Therefore
\[
\tBC \int\limits_{S^\prime} \ \norm{\phi}
\] 
exists and
\[
\tBC \int\limits_{S^\prime} \ \norm{\phi} 
\ \leq \ \tBC \int\limits_S \ \norm{\phi} 
\ < \ 
+\infty,
\]
from which it follows that $\norm{\phi}$ is quasi additive on $S^\prime$.  
Given $\varepsilon > 0$, determine the parameters in the definition of quasi additive in such a way that 
the relevant relations are simultaneously satisfied per $\phi$ on $S$ and per $\norm{\phi}$ on $S^\prime$, 
hence
\begin{align*}
\sum\limits_{I_0} s(I_0,S^\prime) \ 
\normx{\sum\limits_I s(I,I_0) \ \phi(I) - \phi(I_0)} \ 
&\leq \ 
\sum\limits_{I_0} s(I_0,S) \ 
\normx{
\sum\limits_I s(I,I_0) \ 
\phi(I) - \phi(I_0)} 
\\[15pt]
& < \  
\varepsilon
\end{align*}
and
\[
\sum\limits_{I} s(I,S^\prime) \ 
\bigg[
1 - \sum\limits_{I_0} \ 
s(I,I_0) \hsy s(I_0,S^\prime) 
\bigg] 
\norm{\phi(I)} 
\ < \  
\varepsilon.
\]
\\[-.75cm]

\noindent
Therefore $\phi$ is quasi additive on $S^\prime$.
\\[-.25cm]
\end{x}

\begin{x}{\small\bf APPLICATION} \ 
If $\phi:\sI \ra \R^M$ is quasi additive and if 
\[
\tBC \int\limits_A \ \norm{\phi} \ < \ +\infty,
\] 
then $\phi$ is quasi additive on every subset of $A$.
\\[-.25cm]
\end{x}

Here is a summary of the fundamental points of this $\S$.  Work with $\phi$ and $\norm{\phi}$.
\\[-.25cm]

\qquad \textbullet \quad Suppose that $\norm{\phi}$ is quasi subadditive on $S$ and 
\[
\tBC \int\limits_S \ \norm{\phi} \ < \ +\infty.
\] 
Then $\norm{\phi}$ is quasi additive on  $S$ .
\\

\qquad \textbullet \quad Suppose that ${\phi}$ is quasi additive on $S$ 
$-$then  $\norm{\phi}$ is quasi subadditive on  $S$ .
\\

So$:$\quad If $\phi$ is quasi additive on $S$ AND if
\[
\tBC \int\limits_S \ \norm{\phi} \ < \ +\infty,
\] 
then $\norm{\phi}$ is quasi additive on $S$.
\\[-.5cm]

[Note: \  It is not true in general that $\norm{\phi}$ quasi additive implies $\phi$ quasi additive.]
\\[-.25cm]

\begin{x}{\small\bf EXAMPLE} \ 
Take $A = [a,b]$ and let $\sI, \ \sD$, and $\delta$ be as at the beginning.  
Given a continuous curve
\[
C \longleftrightarrow \un{f}:[a,b] \ra \R^M,
\]
define a quasi additive interval function $\phi:\sI \ra \R^M$ by the rule
\begin{align*}
\phi(I) \ 
&= \ 
(\phi_1(I), \ldots, \phi_M(I)) 
\\[8pt]
&= \ 
(f_1(d) - f_1(c), \ldots, f_M(d) - f_M(c)),
\end{align*}
where $I = [c,d] \subset [a,b]$, thus

\[
\norm{\phi(I)} 
\ = \  
\norm{\un{f}(d) - \un{f}(c)},
\]
so if $P \in \sP[a,b]$ corresponds to
\[
D \longleftrightarrow \{[x_{i-1},x_i]: i = 1, \ldots, n\},
\]
then
\[
\sum\limits_{I \in D} \ 
\norm{\phi(I)} 
\ = \ 
\sum\limits_{i=1}^n \ 
\norm{\un{f}(x_i) - \un{f}(x_{i-1})}
\]

$\implies$
\begin{align*}
\tBC \int\limits_A \norm{\phi} \ 
&= \ 
\lim\limits_{\delta(D) \ra 0} \ \sum\limits_{I \in D} \ \norm{\phi(I)}
\\[15pt]
&= \ 
\lim\limits_{\norm{P} \ra 0} \ 
\sum\limits_{i = 1}^n \ \norm{\un{f}(x_i) - \un{f}(x_{i-1})}
\\[15pt]
&= \ 
\ell(C).
\end{align*}
Therefore $C$ is rectifiable iff
\[
\tBC \int\limits_A \norm{\phi} 
\ < \ 
+\infty.
\]
And when this is the case, $\norm{\phi}$ is quasi additive on $A$.
\\[-.5cm]

[Note: \  A priori,
\[
\ell(C) 
\ = \  
\sup\limits_{P \in \sP[a,b]} \ 
\sum\limits_{i=1}^n \ 
\norm{\un{f}(x_i) - \un{f}(x_{i-1})}.
\]
But here, thanks to the continuity of $\un{f}$, the $\sup$ can be replace by $\lim$.]
\\[-.25cm]
\end{x}

\begin{x}{\small\bf EXAMPLE} \ 
Take $A = [a,b]$ and let $\sI$ and $\sD$ be as above.  Suppose that
\[
C \longleftrightarrow  \un{f}:[a,b] \ra \R^M
\]
is a rectifiable curve, potentially discontinuous.
\\

\qquad \textbullet \quad Given $a \leq x_0 < b$, put
\[
s^+(x_0) 
\ = \  
\limsup\limits_{x\downarrow x_0}\  \norm{\un{f}(x) - \un{f}(x_{0})}
\]
and let $s^+(b) = 0.$
\\

\qquad \textbullet \quad Given $a < x_0 \leq b$,  put
\[
s^-(x_0) 
\ = \  
\limsup\limits_{x\uparrow x_0} \norm{\un{f}(x) - \un{f}(x_{0})}
\]
and let $s^-(a) = 0$.  Combine the data and set
\[
s(x) 
\ = \  
s^+(x) + s^-(x) \qquad (a \leq x \leq b).
\]
Then $s(x)$ is zero everywhere save for at most countably many $x$ and 
\[
\sigma 
\ = \  
\sum\limits_x s(x) 
\ \leq \  
\ell(C).
\]
Take $\phi$ as above and define a mesh $\delta$ by the rule
\[
\delta(D) 
\ = \  
\norm{P} + \sigma - \sum\limits_{i=0}^n s(x_i).
\]
One can then show that $\phi$ is quasi additive and
\[
\tBC \int\limits_A \norm{\phi} 
\ = \  
\ell(C).
\]
\\[-1cm]
\end{x}

\begin{x}{\small\bf NOTATION} \ 
Given a quasi additive interval function $\phi$, let
\[
V [\phi,S] 
\ = \  
\sup\limits_{D \in \sD} \  \sum  \ [\hsy \norm{\phi},\hsy S,\hsy D \hsy].
\]
\\[-1.5cm]
\end{x}

\begin{x}{\small\bf \un{N.B.}} \ 
By definition,
\[
\tBC \int\limits_S \ \norm{\phi} 
\ = \  
\lim\limits_{\delta(D) \ra 0} \  \sum  \ [\hsy \norm{\phi},\hsy S,\hsy D \hsy],
\]
so
\[
\tBC \int\limits_S \norm{\phi} 
\ \leq \  
V [\phi,S]
\]
and strict inequality may hold.
\\[-.25cm]
\end{x}

\begin{x}{\small\bf LEMMA} \ 
Given a quasi additive $\phi$ and a subset $S \subset A$, 
suppose that for every $\varepsilon > 0$ and any $D_0 = [I_0]$ 
there exists $\lambda(\varepsilon,S,D_0) > 0$ such that for every system $D = [I]$ with 
$\delta(D) < \lambda(\varepsilon,S,D_0)$ the relation
\[
\sum\limits_{I_0} s(I_0,S) \ 
\bigg[
\sum\limits_I s(I,I_0) \ 
\norm{\phi(I)} - \norm{\phi(I_0)}
\bigg]^- 
\ < \ \varepsilon
\]
obtains $-$then
\[
\tBC  \int\limits_S \ \norm{\phi} 
\ = \ 
V [\phi,S] .
\]
\end{x}

\chapter{
$\boldsymbol{\S}$\textbf{10}.\quad LINE INTEGRALS (bis)}
\setlength\parindent{2em}
\setcounter{theoremn}{0}
\renewcommand{\thepage}{\S10-\arabic{page}}

\qquad 
Throughout this $\S$, the situation will be absolute, where $A = [a,b]$ and $\sI, \  \sD$, and $\delta$ have their usual connotations.
\\

If
\[
C \longleftrightarrow \un{f}:[a,b] \ra \R^M
\]
is a curve, continuous and rectifiable, then
\[
\tBC \ \int\limits_A \ \norm{\phi} \ = \ \ell(C).
\]
And if $F$ is a parametric integrand, then
\[
\int\limits_C F 
\ = \ 
\lim\limits_{\norm{P} \ra 0} \ 
\sum\limits_{i=1}^n \ 
F(\un{f}(\xi_i),\un{f}(x_i) - \un{f}(x_{i-1}))
\]
exists, the result being independent of the $\xi_i$.
\\

\begin{x}{\small\bf \un{N.B.}} \ 
Recall the procedure:  Introduce the integral
\[
I(C) 
\ = \ 
\int\limits_a^b \ 
F(\un{f}(x), \un{\theta}(x)) \ \td\mu_C
\]
and prove that
\[
\lim\limits_{\norm{P} \ra 0} \ 
\sum\limits_{i=1}^n \ 
F(\un{f}(\xi_i),\un{f}(x_i) - \un{f}(x_{i-1}))
\]
exists and equals $I(C)$, the result being denoted by the symbol
\[
\int\limits_C \ F
\]
and called the line integral of $F$ along $C$.
\\[-.25cm]
\end{x}

There is another approach to all of this which does not use measure theory.  
Thus define an interval function $\Phi: \sI \ra \R$ by the prescription
\[
\Phi(I;\xi) \ = \ F(\un{f}(\xi), \phi(I)), 
\]
where $\xi \in I$ is arbitrary.
\\[-.5cm]

[Note: \quad By definition
\[
\phi(I) = (\phi_1(I), \ldots, \phi_M(I)) \ = \ (f_1(d)-f_1(c), \ldots, f_M(d)-f_M(c)),
\]
$I$ being $[c,d] \subset [a,b]$.  
Moreover, $\phi$ is quasi additive.]
\\[-.25cm]

\begin{x}{\small\bf THEOREM} \ 
$\Phi$ is quasi additive. 
\\

Admit the contention $-$then
\[
\lim\limits_{\delta(D)\ra 0} \ \sum\limits_{I \in D} \ \Phi(I;\xi)  
\ = \ 
\lim\limits_{\norm{P} \ra 0} \ 
\sum\limits_{i=1}^n \ 
F(\un{f}(\xi_i),\un{f}(x_i) - \un{f}(x_{i-1}))
\]
exists, call it
\[
(\xi) \ \int\limits_C F.
\]
\\[-1cm]
\end{x}

\begin{x}{\small\bf \un{N.B.}} \ 
Needless to say, it turns out that
\[
(\xi) \ \int\limits_C F
\]
is independent of the $\xi$ $($this follows by a standard ``$\varepsilon/3$'' argument$)$ $($details at the end$)$.
\\[-.25cm]

[Note: \ This is one advantage of the approach via $I(C)$ in that independence is manifest.]
\\[-.25cm]
\end{x}

To simplify matters, it will be best to generalize matters.
\\[-.25cm]

Assume from the outset that $\phi:\sI \ra \R^M$ is now an arbitrary interval function which is quasi additive with 
\[
\tBC\ \int\limits_A \norm{\phi}  
\ < \ +\infty,
\]
hence that $\norm{\phi}$ is also quasi additive as well.
\\[-.25cm]

Introduce another interval function $\zeta:\sI \ra \R^N$ and expand the definition of parametric integrand so that
\[
F:X \times \R^M \ra \R,
\]
where $X \subset \R^N$ is compact and $\zeta(I) \subset X$.
\\[-.25cm]

\begin{x}{\small\bf EXAMPLE} \ 
To recover the earlier set up, take $N = M$, keep $\phi:\sI \ra \R^M$, let $\omega:\sI \ra [a,b]$ be a choice function, i.e., suppose that $\omega(I) \in I \subset [a,b]$, let $\zeta(I) = \un{f}(\omega(I))$, and take $X = [C] \subset \R^M$.
\\[-.25cm]
\end{x}

\begin{x}{\small\bf CONDITION} \ 
$(\zeta)$ \quad $\forall \ \varepsilon > 0$, $\exists \ t(\varepsilon) > 0$ such that if $D_0 = [I_0]$ is any system subject to $\delta(D_0) < t(\varepsilon)$ there also exists $T(\varepsilon, D_0)$, such that for any system $D = [I]$ with $\delta(D) < T(\varepsilon, D_0)$ the relation
\[
\max\limits_{I_0} \  \max\limits_{I \subset I_0}\  \norm{\zeta(I) - \zeta(I_0)} 
\ < \ \varepsilon
\]
obtains.
\\[-.25cm]
\end{x}

\begin{x}{\small\bf \un{N.B.}} \ 
Owing to the uniform continuity of $\un{f}$, this condition is automatic in the special case supra.
\\[-.25cm]
\end{x}

\begin{x}{\small\bf THEOREM} \ 
Let $F$ be a parametric integrand, form the interval function $\Phi:\sI \ra \R$ defined by the prescription
\[
\Phi(I) \ = \ F(\zeta(I),\phi(I)),
\]
and impose condition $(\zeta)$ $-$then $\Phi$ is quasi additive.
\\[-.25cm]

The proof will emerge from the discussion below but there are some preliminaries that have to be dealt with first.
\\[-.5cm]
\end{x}

Start by writing down simultaneously ($\tqa_1$-$A$) and ($\tqa_2$-$A$) for 
$\phi$ and $\norm{\phi}$ (both are quasi additive), $\bar{\varepsilon}$ to be determined.

\[
\begin{cases}
\ds\ \sum\limits_{I_0} \ 
\normxx{
\sum\limits_{I \subset I_0} \ \phi(I) - \phi(I_0)
}
\ < \ \bar{\varepsilon}
\\[26pt]
\ds\ \sum\limits_{I \not\subset I_0} \ \norm{\phi(I)} \ < \ \bar{\varepsilon}
\end{cases}
\quad
\begin{cases}
\ds\ \sum\limits_{I_0} \ 
\bigg|
\sum\limits_{I \subset I_0} \ \norm{\phi(I)} - \norm{\phi(I_0)}
\bigg|
\ < \ \bar{\varepsilon}
\\[26pt]
\ds\ \sum\limits_{I \not\subset I_0} \ \abs{\norm{\phi(I)}} \ < \ \bar{\varepsilon}
\end{cases}
\]
for $\delta(D_0) < \eta(\bar{\varepsilon})$ and $\delta(D) < \lambda(\bar{\varepsilon}, D_0)$ and in addition
\[
\bigg| 
\sum\limits_{I \in D} \ \norm{\phi(I)} - \tBC \ \int\limits_A \ \norm{\phi(I)} 
\bigg|
\ < \ \bar{\varepsilon}
\]
for $\delta(D) < \sigma(\bar{\varepsilon})$.\\

Fix $\varepsilon > 0$.  Put
\[
V 
\ = \  
\tBC \ \int\limits_A \norm{\phi} \qquad (< +\infty).
\]
\\[-.5cm]

\qquad \textbullet \quad 
$(F)$ \quad $X \times \tU(M)$ is a compact set on which $F$ is bounded:
\[
\abs{F(\un{x},\un{t})} \ \leq \ C \qquad(\un{x} \in X,\hsy \un{t} \in \tU(M))
\]
and uniformly continuous: $\exists \ \gamma$ such that

\[
\begin{cases}
\norm{\un{x} - \un{x}^\prime}
\\[8pt]
\norm{\un{t} - \un{t}^\prime}
\end{cases}
< \ \gamma 
\ \implies \ 
\abs{F(\un{x},\un{t}) - F(\un{x}^\prime,\un{t}^\prime)} \ < \  \frac{\varepsilon}{3(V + \varepsilon)}.
\]
\\[-.25cm]

\qquad \textbullet \quad 
$(\alpha)$

\[
\alpha(I_0) \ = \ \frac{\phi(I_0)}{\norm{\phi(I_0)}} \quad \text{if } \phi(I_0) \neq 0
\]

\noindent
but 0 otherwise and

\[
\alpha(I) \ = \ \frac{\phi(I)}{\norm{\phi(I)}} \quad \text{if } \phi(I) \neq 0
\]
but 0 otherwise.
\\

\begin{x}{\small\bf NOTATION} \ 
Denote by
\[
\sum {}_{\gamma^+}^{(I_0)}
\]
the sum over the $I \subset I_0$ for which
\[
\norm{\alpha(I_0) - \alpha(I)} \ \geq \ \gamma
\]
and denote by
\[
\sum {}_{\gamma^-}^{(I_0)}
\]
the sum over the $I \subset I_0$ for which
\[
\norm{\alpha(I_0) - \alpha(I)} \ < \ \gamma.
\]

Therefore
\[
\sum\limits_{I \subset I_0} \ = \ 
\sum {}_{\gamma^+}^{(I_0)} 
\ + \ 
\sum {}_{\gamma^-}^{(I_0)}.
\]
\\[-1cm]
\end{x}

\begin{x}{\small\bf LEMMA} \ 
\\[-1cm]

\begin{align*}
\frac{\gamma^2}{2} \ 
\sum\limits_{I_0} \ 
\sum {}_{\gamma^+}^{(I_0)} 
\ \norm{\phi(I)} 
\ \leq \  \ 
&
\sum\limits_{I_0} \ 
\Big\Vert
\sum\limits_{I \subset I_0} \ 
\phi(I) - \phi(I_0)
\Big\Vert \ 
 + \ 
 \sum\limits_{I_0} \ 
 \Big|
 \sum\limits_{I \subset I_0} \ 
 \norm{\phi(I)} - \norm{\phi(I_0)} 
 \Big|.
\end{align*}
\\[-.75cm]

PROOF \ 
The inequality
\[
\norm{\alpha(I_0) - \alpha(I)} \ \geq \ \gamma
\]
implies that
\begin{align*}
\gamma^2 \ 
&\leq \ 
\norm{\alpha(I_0) - \alpha(I)}^2
\\[11pt]
&=\ 
 (\alpha(I_0) - \alpha(I)) \cdot (\alpha(I_0) - \alpha(I))
\\[11pt]
&= \ 
\norm{\alpha(I_0)}^2 - 2\alpha(I_0) \cdot \alpha(I) + \norm{\alpha(I)}^2
\\[11pt]
&= 2 - 2\alpha(I_0) \cdot \alpha(I),
\end{align*}
so
\[
\frac{\gamma^2}{2} \ \leq \ 1 - \alpha(I_0) \cdot \alpha(I)
\]

\qquad
$\implies$
\[
\frac{\gamma^2}{2} \norm{\phi(I)} \ \leq \ \norm{\phi(I)} - \alpha(I_0) \cdot \phi(I).
\]
But for any $I$, 
\[
0 \ \leq \ \norm{\phi(I)} - \alpha(I_0) \cdot \phi(I).
\]
Proof: \ In fact,
\[
\norm{\phi(I)} - \frac{\phi(I_0) \cdot \phi(I)}{\norm{\phi(I_0)}} 
\ = \ 
\frac{1}{\norm{\phi(I_0)}} 
\left[
\norm{\phi(I)} \norm{\phi(I_0)} - \phi(I_0) \cdot \phi(I)
\right]
.
\]

\noindent
Now quote Schwarz's inequality.  
Thus we may write
\begin{align*}
\frac{\gamma^2}{2} \ 
\sum {}_{\gamma^+}^{(I_0)} \norm{\phi(I)} \ 
&\leq \ 
\sum {}_{\gamma^+}^{(I_0)} \ 
(\norm{\phi(I) - \alpha(I_0) \cdot \phi(I))}
\\[15pt]
&\leq \ 
\sum\limits_{I \subset I_0} \ 
(\norm{\phi(I) - \alpha(I_0) \cdot \phi(I))} 
\\[15pt]
&= \ 
\sum\limits_{I \subset I_0} \ 
\norm{\phi(I)} - \norm{\phi(I_0)} 
+ \alpha(I_0) \cdot \bigl(\phi(I_0) - \sum\limits_{I \subset I_0} \ \phi(I)\bigr)
\\[15pt]
&\leq \ 
\bigg|\sum\limits_{I \subset I_0} \ 
(\norm{\phi(I)} - \norm{\phi(I_0)}\bigg| 
+ 
\normx
{\phi(I_0) - \sum\limits_{I \subset I_0} \ \phi(I)}
\qquad (\text{Schwarz}).
\end{align*}
To finish, sum over $I_0$.
\end{x}

\qquad \textbullet \quad 
$(D_0)$ \quad Assume
\[
\qquad \qquad \delta(D_0) 
\ < \ 
\min\{t(\gamma), \hsy\eta(\varepsilon), \hsy\eta(\varepsilon\gamma^2)\}.
\]

\qquad \textbullet \quad 
$(D)$ \quad Assume
\[
\qquad \qquad  \delta(D) 
\ < \ 
\min\{\sigma(\varepsilon), \hsy\lambda(\varepsilon,D_0), \hsy\lambda(\varepsilon\gamma^2,D_0), \hsy T(\gamma,D_0)\}.
\]

\qquad \textbullet \quad 
$(\bar{\varepsilon})$ \quad Assume
\[
\qquad \qquad  \bar{\varepsilon} 
\ < \ 
\min\left\{\gamma,\hsy\frac{\varepsilon}{3C},\hsy \frac{\varepsilon\gamma^2}{24 \hsy C}\right\}.
\]
Then

\begingroup
\allowdisplaybreaks
\begin{align*}
\sum\limits_{I_0} \ 
\Big|
\sum\limits_{I \subset I_0} \ 
&\Phi(I) - \Phi(I_0)
\Big| \  
\\[15pt]
&= \ 
\sum\limits_{I_0} \ 
\Big|
\sum\limits_{I \subset I_0} \ 
F(\zeta(I),\phi(I)) 
-  
\sum\limits_{I\subset I_0} \ 
F(\zeta(I_0),\alpha(I_0))\norm{\phi(I)}
\\[18pt]
&\qquad\qquad +  
\sum\limits_{I \subset I_0} \ 
F(\zeta(I_0),\alpha(I_0))\norm{\phi(I)} - F(\zeta(I_0),\alpha(I_0))\norm{\phi(I_0)} 
\Big|
\\[18pt]
&= \ 
\sum\limits_{I_0} \ 
\Big|
\sum\limits_{I \subset I_0} \ 
F(\zeta(I),\alpha(I)) -   F(\zeta(I_0),\alpha(I_0)))
\norm{\phi(I)}
\\[18pt]
&\qquad\qquad 
+  
\sum\limits_{I \subset I_0} \ 
F(\zeta(I_0),\alpha(I_0)) \hsy
(\norm{\phi(I)} -\norm{\phi(I_0)}) 
\Big|
\\[18pt]
&\leq  \ 
\sum\limits_{I_0} \Big|F(\zeta(I_0),\alpha(I_0))\Big| \quad \Big| \sum\limits_{I \subset I_0} \norm{\phi(I)} - \norm{\phi(I_0)}  \Big|
\\[18pt]
&\qquad\qquad 
+  
\sum\limits_{I_0}\ 
\sum\limits_{I \subset I_0} \ 
|F(\zeta(I),\alpha(I)) - F(\zeta(I_0),\alpha(I_0))| \ \norm{\phi(I)} 
\\[18pt]
&= \  
\sum\limits_{I_0} \ 
\Big|
F(\zeta(I_0),\alpha(I_0))
\Big| 
\quad 
\Big| 
\sum\limits_{I \subset I_0} \
\norm{\phi(I)} - \norm{\phi(I_0)}  
\Big|
\\[18pt]
&\qquad + \sum\limits_{I_0} \ 
\Big(
\sum {}_{\gamma^-}^{(I_0)} + \sum {}_{\gamma^+}^{(I_0)} 
\Big) 
\big| 
F(\zeta(I_0),\alpha(I_0)) - F(\zeta(I),\alpha(I)) 
\big| 
\norm{\phi(I)}.
\end{align*}

\un{First}:

\begin{align*}
\sum\limits_{I_0} \Big|F(\zeta(I_0),\alpha(I_0))
\Big| 
\ \hsx
\Big| 
\sum\limits_{I \subset I_0} \norm{\phi(I)} - \norm{\phi(I_0)}  
\Big| \ 
&\leq \ 
C \hsx \sum\limits_{I_0} \ 
\Big| 
\sum\limits_{I \subset I_0} \ 
\norm{\phi(I)} - \norm{\phi(I_0)}  
\Big|
\\[15pt]
&< \ 
C \bar{\varepsilon} 
\\[15pt]
&< \ 
C \frac{\varepsilon}{3 C}
\\[15pt]
&= \ 
\frac{\varepsilon}{3}.
\end{align*}
\\[-.75cm]

\un{Second}:  Consider

\[
\sum\limits_{I_0} \ 
\sum {}_{\gamma^-}^{(I_0)} \ 
\big| F(\zeta(I_0),\alpha(I_0)) - F(\zeta(I),\alpha(I)) \big| \norm{\phi(I)}.
\]
Here

\[
\begin{cases}
\norm{\alpha(I_0)} = 1, 
\\[8pt] 
\norm{\alpha(I)} = 1, 
\\[8pt] 
\norm{\alpha(I_0) - \alpha(I)} < \gamma, 
\\[8pt] 
\norm{\zeta(I_0) - \zeta(I)} < \gamma
\end{cases}
\implies
\big| F(\zeta(I_0),\alpha(I_0)) - F(\zeta(I),\alpha(I)) \big| 
\ < \ 
\frac{\varepsilon}{3(V+\varepsilon)}.
\]
\\[-.5cm]

\noindent
The entity in question is thus majorized by
\begin{align*}
\frac{\varepsilon}{3(V+\varepsilon)} \ 
\sum\limits_{I_0} \ 
\sum {}_{\gamma^-}^{(I_0)} \norm{\phi(I)} \ 
&\leq \ 
\frac{\varepsilon}{3(V+\varepsilon)} \ 
\sum\limits_{I \in D} \ 
\norm{\phi(I)}
\\[15pt]
&\leq \ 
\frac{\varepsilon}{3(V+\varepsilon)} (V + \varepsilon)
\\[15pt]
&= \ 
\frac{\varepsilon}{3}.
\end{align*}

\un{Third}:

\begin{align*}
\sum\limits_{I_0} \ 
\sum {}_{\gamma^+}^{(I_0)}  \ 
\big| 
&
F(\zeta(I_0),\alpha(I_0)) - F(\zeta(I),\alpha(I)) 
\big| 
\norm{\phi(I)} 
\\[15pt]
&\leq \ 
2 \hsy C \  \sum\limits_{I_0} \ 
\sum {}_{\gamma^+}^{(I_0)} \ 
\norm{\phi(I)}
\\[15pt]
&\leq \ 
\frac{4 \hsy C}{\gamma^2} \ 
\Big[  
\sum\limits_{I_0}  \ 
\sum\limits_{I \subset I_0} \ 
\norm{\phi(I) - \phi(I_0)} 
\ + \ 
\sum\limits_{I_0} \ 
\Big| 
\sum\limits_{I \subset I_0} \ 
\norm{\phi(I)} - \norm{\phi(I_0)}  
\Big| \hsx
\Big]
\\[15pt]
&\leq \ 
\frac{4 \hsy C}{\gamma^2} (\bar{\varepsilon} + \bar{\varepsilon})
\\[15pt]
&= \ 
\frac{8 \hsy C}{\gamma^2} \bar{\varepsilon}
\\[15pt]
&< \ 
\frac{8 \hsy C}{\gamma^2} \cdot \frac{\varepsilon\gamma^2}{24 \hsy C}
\\[15pt]
&= \  
\frac{\varepsilon}{3}.
\end{align*}
\\[-.75cm]

\noindent
In total then:
\begin{align*}
\sum\limits_{I_0} 
\Big| 
\sum\limits_{I \not\subset I_0} \ 
\Phi(I) - \Phi(I_0) 
\Big| \ 
&< \ \frac{\varepsilon}{3} \hsx + \hsx \frac{\varepsilon}{3} \hsx + \hsx \frac{\varepsilon}{3}
\\[11pt]
&=\  
\varepsilon.
\end{align*}
And finally
\begin{align*}
\sum\limits_{I \not\subset I_0} \ 
\abs{\Phi(I)} \ 
&= \ \sum\limits_{I \not\subset I_0} \ 
\abs{F(\zeta(I),\phi(I))}
\\[15pt]
&= \ 
\sum\limits_{I \not\subset I_0} \ 
\abs{F(\zeta(I),\alpha(I))} \norm{\phi(I)}
\\[15pt]
&\leq \ 
C \hsx \sum\limits_{I \not\subset I_0} \ 
\norm{\phi(I)}
\\[15pt]
&< \ 
C \hsx \bar{\varepsilon} 
\\[15pt]
&< \ 
C \hsx \frac{\varepsilon}{3 \hsy C} 
\\[15pt]
&= \  \frac{\varepsilon}{3} 
\\[15pt]
&< \ 
\varepsilon.
\end{align*}
\endgroup

Therefore $\Phi$ is quasi additive.  
And since the conditions on $F$ carry over to $\abs{F}$, it follows that $\norm{\Phi}$ is also quasi additive, hence
\[
\tBC \ \int\limits_A \ \norm{\Phi}
\]
exists and is finite.
\\[-.25cm]

To tie up one loose end, return to the beginning and consider the line integrals
\[
(\xi) \  \int\limits_C \ F, \qquad (\xi^\prime) \  \int\limits_C \ F,
\]
the claim being that they are equal.  That this is so can be seen by writing
\begin{align*}
\Big|(\xi)  \ \int\limits_C \ F 
\ \  - \  \  
&
(\xi^\prime) \  \int\limits_C \ F \ \Big| \ 
\\[15pt]
&= \ 
\Big|(\xi)  \int\limits_C \ F - \sum\limits_{i=1}^n F(\un{f}(\xi_i),\un{f}(x_i) - \un{f}(x_{i-1}))
\\[15pt]
&\quad\qquad + \sum\limits_{i=1}^n F(\un{f}(\xi_i),\un{f}(x_i) - \un{f}(x_{i-1}))
\\[15pt]
&\quad\qquad -  \sum\limits_{i=1}^n F(\un{f}(\xi_i^\prime),\un{f}(x_i) - \un{f}(x_{i-1}))
\\[15pt]
&\quad\qquad 
+  
\sum\limits_{i=1}^n \ 
F(\un{f}(\xi_i^\prime),\un{f}(x_i) - \un{f}(x_{i-1})) \ - \ (\xi^\prime)  \int\limits_C \ F \ \Big|
\end{align*}
and proceed from here in the obvious way.
\\[-.25cm]

\begin{x}{\small\bf EXAMPLE} \ 
Take $N = 1$, $M = 1$ and define an interval function $\abs{\ \cdot\ }: \sI \ra \R$ by sending $I$ to its length $\abs{I}$.  
Fix a choice function $\omega : \sI \ra [a,b]$.  
Consider a curve
\[
C \longleftrightarrow f:[a,b] \ra \R.
\]
Assume:\quad $f$ is continuous and of bounded variation, thus
\[
\ell(C) \ = \  T_f[a,b] \ < \ +\infty.
\]
Work with the parametric integrand $F(x,t) = xt$ $-$then the data
\begin{align*}
I 
&\ra 
F(\zeta(I),\abs{I})
\\[11pt]
&=\  
F(f(\omega(I)),\abs{I})
\\[11pt]
&=\  
f(\omega(I)) \abs{I}
\end{align*}
leads to sums of the form
\[
\sum\limits_{i=1}^n \ f(\xi_i)(x_i - x_{i-1}),
\]
hence to
\[
\int\limits_C F 
\ = \  
\int\limits_a^b f,
\]
the Riemann integral of $f$.
\end{x}

\chapter{
$\boldsymbol{\S}$\textbf{11}.\quad EXAMPLES}
\setlength\parindent{2em}
\setcounter{theoremn}{0}
\renewcommand{\thepage}{\S11-\arabic{page}}

\begin{x}{\small\bf EXAMPLE} \
Suppose that $A = [0,1]$ with the usual topology, 
$\{I\}$ the class of all intervals $I = [a,b] \subset A$, $a < b$, 
and define $\abs{I} = b - a$.  
Let $\sD$ be the family of all finite systems $D = [I]$ of nonoverlapping $I \in \{I\}$.  
If we take 
\[
\delta (D) 
\ = \ 
\left(
1 - \sum \ \abs{I}
\right) 
+ \max \abs{I},
\]
then
\[
0 
\ < \ 
\delta (D) < 2,
\]
and the Assumption \S9, \#4, is trivially satisfied.

Take $k = 1$ and $\phi(I) = b - a > 0$.  
Given $\varepsilon > 0$ arbitrary take $\eta = \eta (\varepsilon) = \ds\frac{\varepsilon}{2}$ and let 
\[
D_0 
\ = \ 
(I_j, \ j = 1, \ldots,  N)
\]
be any system $D_0 \in \sD$ with $\delta(D_0) < \eta$.  
Take 
\[
\lambda
\ = \ 
\lambda (\varepsilon, D_0) 
\ = \ 
\frac{\varepsilon}{4 N}.
\]
If $D = [J]$ is any system with $\delta (D) < \lambda$, 
and we denote by $\sum {}^\prime$ any sum ranging over all $J \in D$ with $J \not\subset I$ and for any $I \in D_0$, 
we have
\begin{align*}
\sum {}^\prime \ 
\phi(J) \ 
&< \ 
\left(
1 - \sum \abs{I_j}
\right)
+ 2\hsy  N(\varepsilon/4 N) 
\\[15pt]
&< \ 
\frac{\varepsilon}{2} + \frac{\varepsilon}{2} 
\\[15pt]
&=\ 
\varepsilon, 
\end{align*}
and $(\phi_2)$ holds.
Also, if $\sum {}^{(j)}$ denotes a sum ranging over all $J \in D$ with $J \subset I_j$, we have 
\[
\phi (I_j) 
- 
\sum {}^{(j)} \ 
\phi (J) 
\ > \ 0,
\]
and
\begin{align*}
 \sum\limits_j \ 
\left(
 \phi (I_j) 
 - 
 \sum {}^{(j)} \  
 \phi (J) 
 \right) \ 
&< \ 
\left(
1 - \sum \abs{J}
\right)
+ 2 \hsy N(\varepsilon/4 N) 
\\[15pt]
&< \ 
\frac{\varepsilon}{2} + \frac{\varepsilon}{2} 
\\[15pt]
&=\ 
\varepsilon, 
\end{align*}
and $(\phi_1)$ holds.  
Thus, $\phi$ is quasi additive with $\sB = 1$.
\\[-.25cm]
\end{x}


\begin{x}{\small\bf EXAMPLE} \ 
Suppose that $A = [0,1]$ with the usual topology, 
$\{I\}$ the class of all intervals $I = [a,b] \subset A$, $a < b$, 
and define $\abs{I} = b - a$.  
Let $\sD$ be the family of all finite systems $D = [I]$ of nonoverlapping $I \in \{I\}$.  
Take 
\[
\delta (D) 
\ = \ 
\begin{cases}
\left(
1 - \sum \ \abs{I}
\right) 
+ \max \abs{I}
\hspace{1.25cm} \text{if for no $I \in D$ we have $a = 0$}
\\[8pt]
\left(
1 - \sum \ \abs{I}
\right) 
+ \max \abs{I} + 1
\hspace{0.5cm} \text{if for one $I \in D$ we have $I = [0,b]$}
\end{cases}
.
\]
If $\eta < 1$, and $\delta (D) < 1$, $D = [I]$, then there is in $D$ no interval $I = [0,b]$.  
The assumptions \S9, \#4 , are trivially satisfied.  
By the same reasoning as in the previous example with 
$\ds\eta = \min \hsx \left[\frac{\varepsilon}{2}, 1\right]$ 
it follows that $\phi$ is quasi additive with $V = 1$.
\\
\end{x}


\begin{x}{\small\bf EXAMPLE} \ 
Suppose that $A = [0,1]$ with the usual topology, 
$\{I\}$ the class of all intervals $I = [a,b] \subset A$, $a < b$, 
and define $\abs{I} = b - a$.  
Let $\sD$ be the family of all finite systems $D = [I]$ of nonoverlapping $I \in \{I\}$.  
Take 
\[
\delta (D) 
\ = \ 
\begin{cases}
\left(
1 - \sum \ \abs{I}
\right) 
+ \max \abs{I} + 1
\hspace{.6cm} \text{if for no $I \in D$ we have $a = 0$}
\\[8pt]
\left(
1 - \sum \ \abs{I}
\right) 
+ \max \abs{I}
\hspace{1.25cm} \text{if for one $I \in D$ we have $I = [0,b]$}
\end{cases}
.
\]
If $\eta < 1$, and $\delta (D) < 1$, $D = [I]$, then there is in $D$ one interval $I = [0,b]$.  
The assumptions \S9, \#4 , are trivially satisfied.  
By the same reasoning as in \#2 we prove that 
$\phi$ is quasi additive with $V = 1$.
\\
\end{x}


\begin{x}{\small\bf EXAMPLE} \
Suppose that $A = [0,1]$ with the usual topology, 
$\{I\}$ the class of all intervals $I = [a,b] \subset A$, $a < b$, 
and define $\abs{I} = b - a$.  
Let $\sD$ be the family of all finite systems $D = [I]$ of nonoverlapping $I \in \{I\}$.  
For every real $m$, $0 \leq m \leq 1$, take 
\[
\sigma (m) 
\ = \ 
\begin{cases}
\ 0 \hspace{0.6cm} \text{if $m$ is irrational}
\\
\ \ds\frac{1}{p} \hspace{0.5cm} \text{if} \ m = \frac{p}{q}, \  p, q \in \Z, \ p \geq 1 \hsx \text{minimum}, \ 0 \leq q \leq p.
\end{cases}
\]
Take $k = 1$ and 
\[
\phi(I) 
\ = \ 
I
\ = \ b - a 
\quad \text{for} \ I = [a,b] \subset A.
\]
For $D = [I_i : j = 1, \ldots, N]$, $I_j = [a_j, b_j]$, take
\[
\delta (D) 
\ = \ 
\left(
1 - \sum \ \abs{I_j}
\right) 
+ \max \abs{I_j} 
+ 
\sum \ \sigma (a_j)
+ 
\sum \ \sigma (b_j)
\]
For $D = [I_1]$, $I_1 = [0,1]$, we have 
\[
\delta (D) 
\ = \ 
0 + 1 + 1 + 1 
\ = \ 
3;
\]
for $D = [I_1, I_2]$, $I_1 = [0,1/2]$,  $I_2 = [1/2,1]$, we have 
\[
\delta (D) 
\ = \ 
1/2 + (1 + 1/2) + (1/2 + 1) 
\ = \ 
3 + 1/2.
\]
Thus $\delta (D)$ does not necessarily decrease by refinement.
The Assumption \S9, \#4 , is trivially satisfied.  
Given $\varepsilon > 0$ let $D = [I_j]$ with $I_j = [a_j, b_j]$, all $a_j$, $b_j$ irrational with 
\[
1 - 
\sum\ 
(b_j - a_j) 
\ < \ 
\frac{\varepsilon}{2},
\quad 
(b_j - a_j) 
\ < \ 
\frac{\varepsilon}{2}.
\]
Then 
\[
\delta (D) < \varepsilon.
\]
We may also choose an integer 
$
M : \ 
\begin{cases}
\ M \geq 3, 
\\[4pt]
\ M > \ds\frac{4}{\varepsilon}
\end{cases}
, M 
$
prime.  
Suppose that $a_j$, $b_j$ are all of the form 
$\ds\frac{q}{M^2}$, $q = 1, 2, \ldots, M^2 - 1$, and $N = M$.  
Namely, we may take
\allowdisplaybreaks
\begin{align*}
a_1 
&= \frac{(M - 1)}{M} + \frac{1}{M^2},
\\[11pt]
a_2 
&= b_1 = \frac{1}{M} + \frac{1}{M^2},
\\[11pt]
a_3 
&= b_2 = \frac{2}{M} + \frac{1}{M^2},
\\[11pt]
\hspace{1cm} \vdots
\\[11pt]
a_N 
&= b_{N-1} = \frac{M - 1}{M} + \frac{1}{M^2},
\\[11pt]
b_N 
&= 1 - \frac{1}{M^2}.
\end{align*}
Then
\begin{align*}
\delta (D) 
&= \ 
\frac{2}{M^2}
\hsx + \hsx 
\frac{1}{M}
\hsx + \hsx 
M \hsy \frac{2}{M^2}
\\[11pt]
&<\ 
\frac{4}{M}
\\[11pt]
&<\ 
\varepsilon.
\end{align*}
This remark shows that the second part of the Assumption \#4, \S9, holds.
The proof that 
($\text{qa}_1$ -A) and ($\text{qa}_2$ -A), (cf. \#14, \S9), 
hold is the same as for \#1.  
Thus $\phi$ is quasi additive and $\sB = 1$.
\\
\end{x}


\begin{x}{\small\bf EXAMPLE} \ 
Suppose that $A = [0,1]$ with the usual topology, 
$\{I\}$ the class of all intervals $I = [a,b] \subset A$, $a < b$, 
and define $\abs{I} = b - a$.  
Let $\sD$ be the family of all finite systems $D = [I]$ of nonoverlapping $I \in \{I\}$.  
Suppose that $k = 1$ and 
\[
\phi (I) \ = \ 
\begin{cases}
\ \abs{I} = b -a \hspace{0.5cm} \text{ if both $a$ and $b$ are irrational}
\\[4pt]
\ \  a - b \hspace{1.5cm} \text{otherwise}
\end{cases}
.
\]
For 
$D = [I] \in \sD$, $I = [a,b],$
take
\[
\delta (D) 
\ = \ 
\begin{cases}
\left(
1 - 
\sum\ 
\abs{I}
\right)
 + \max \abs{I}
\hspace{0.5cm} \text{ if both $a$ and $b$ are irrational}
 \\[11pt]
 \ 1 \hspace{4.25cm} \text{otherwise}
 \end{cases}
 .
\]
Obviously the Assumptions \S9, \#4 are satisified.  
If $D = [I]$, $\delta (D) < 1$, then all $a$ and $b$ are irrational.  
Obviously ($\text{qa}_1$ -A) and ($\text{qa}_2$ -A), (cf. \#14 \S9), are satisfied.  
$\phi (I)$ is quasi additive, and $\sB = 1$.  
\\
\end{x}


\begin{x}{\small\bf EXAMPLE} \ 
Suppose that $A = [0,1]$ with the usual topology, 
$\{I\}$ the class of all intervals $I = [a,b] \subset A$, $a < b$, 
and define $\abs{I} = b - a$.  
Let $\sD$ be the family of all finite systems $D = [I]$ of nonoverlapping $I \in \{I\}$.  
Suppose that $k = 1$ and 
\[
\phi (I) \ = \ 
\begin{cases}
\ \abs{I} = b -a \hspace{0.5cm} \text{ if both $a$ and $b$ are irrational}
\\[4pt]
\ \  a - b \hspace{1.5cm} \text{otherwise}
\end{cases}
.
\]
For 
$D = [I] \in \sD$, $I = [a,b],$
take
\[
\delta (D) 
\ = \ 
\begin{cases}
\left(
1 - 
\sum\ 
\abs{I}
\right)
 + \max \abs{I}
\hspace{0.5cm} \text{ if both $a$ and $b$ are rational}
 \\[11pt]
 \ 1 \hspace{4.25cm} \text{otherwise}
 \end{cases}
 .
\]
Obviously the Assumptions \S9, \#4 are satisified.  
If $D = [I]$, $\delta (D) < 1$, then all $a$ and $b$ are irrational.  
Obviously ($\text{qa}_1$ -A) and ($\text{qa}_2$ -A), (cf. \#14 \S9), are satisfied.  
$\phi (I)$ is quasi additive, and $\sB = -1$.  
\\
\end{x}


\begin{x}{\small\bf EXAMPLE} \ (Jordan length for continuous curves) \ 
\\[-.5cm]

Let $C : x = x(t)$, $0 \leq t \leq 1$, $x = (x_1, \ldots, x_k)$, 
be any real continuous vector function (a continuous parameteric curve).   
Take $A = [0,1]$, $\{I\}$ the class of all intervals $I = [a,b] \subset A$, 
$\sD$ the class of all finite subdivisions 
$D = [I_j, \hsy j = 1, \ldots , N]$, $I_j = [a_{j - 1}, a_j]$, $0 = a_0 < a_1 < \cdots < a_N = 1$, 
and suppose that 
\[
\delta (D)
\ = \  
\max \abs{I_j}
\ = \  
\max \{a_j - a_{j - 1}\}.  
\]
Then
\[
0 
\ < \ 
\delta (D) 
\ \leq \ 
1
\]
and the assumptions \S9, \#4 , are trivially satisfied.
Suppose that 
\[
\begin{cases}
\ \phi (I) = (\phi_1, \ldots, \phi_k)
\\
\quad \text{where} 
\\
\ \phi_r = x_r (b) - x_r(a) \quad (r = 1, \ldots, k)
\end{cases}
\forall \ I = [a,b] \subset A.  
\]
Denote by 
\[
\begin{cases}
\ M_r = \max\limits_{0 \leq t \leq 1} \abs{x_t(t)} 
\\[4pt]
\ M = \max\limits_{1 \leq r \leq k} \  M_r
\end{cases}
.
\]
Given $\varepsilon > 0$ take $\eta =  \eta (\varepsilon) = 2$.  
Let 
\[
D_0 
\ = \ 
[I_i : i = 1, \ldots, N]
\]
be any system $D_0 \in \sD$ (with $\delta (D_0) < 2)$.  
Because of the continuity of $x(t)$ there is a 
$\lambda = \lambda (\varepsilon, D_0) > 0$ such that $\ds \abs{\phi (I)} < \frac{\varepsilon}{2 N}$ for all $I = [a,b]$ with $b - a < \lambda$.  
Now suppose $D = \abs{J}$ is any system $D \in \sD$ with $\delta (D) < \lambda$.  
Then the sum 
$\sum{}^\prime \ \norm{\phi(J)}$ 
contains at most $2 (N - 1)$ terms and 
\[
\sum{}^\prime \ \norm{\phi(J)}
\ < \ 
2 (N - 1) \frac{\varepsilon}{2 N}
\ < \ 
\varepsilon.
\]
For each $I_i \in D_0$ the intervals $J \in D$, $J \subset I_i$ leave uncovered in $I$ at most two terminal intervals, 
say $H^\prime$, $H^{\prime\prime}$ (if any).  
Then we have 
\[
\phi (I_i) 
\ = \ 
\sum{}^{(i)} \ \phi(J) 
\hsx + \hsx \phi (H^\prime) 
\hsx + \hsx \phi (H^{\prime\prime}),
\]
\begin{align*}
\normx{\phi (I_i)  - \sum{}^{(i)} \ \phi(J)} \ 
&\leq \ 
\norm{\phi (H^\prime) } + \norm{\phi (H^{\prime\prime})}
\\[11pt]
&<\ 
\frac{\varepsilon}{N},
\end{align*}
and
\begin{align*}
\sum\limits_{i = 1, \ldots, N}
\normx{\phi (I_i)  - \sum\limits_{\substack{J \in D\\ J \subset I_i}}{}^{(i)} \ \phi(J)} \ 
&< \ 
N \hsy \frac{\varepsilon}{N}
\\[11pt]
&= \
\varepsilon.
\end{align*}
Thus $\phi$ is quasi additive, as are $\phi_r$, 
$\norm{\phi}$, 
$\phi_r^+$, $\phi_r^-$ 
(cf. \S9, \#16).  
Obviously 
$V = V(\norm{\phi})$ 
is the Jordan length of the curve $C$, and 
$V_r = V(\abs{\phi_r})$, 
$V_r^+ = V(\phi_r^+)$, 
$V_r^- = V(\phi_r^-)$, 
are the total variation, the positive and negative variations of $x_r(t)$, $0 \leq t \leq 1$.  
And
\[
V < +\infty
\implies
\begin{cases}
\ \text{$C$ is rectifiable}
\\[4pt]
\ V_r, \hsy V_r^+, \hsy V_r^- \hsx < \hsx +\infty 
\\[4pt]
\ \phi,  \hsy \phi_r, \hsy \norm{\phi}, \hsy \abs{\phi}, \hsy \phi_r^+, \hsy \phi_r^- 
\quad \text{are all quasi additive}
\end{cases}
.
\]
\\[-1cm]
\end{x}


\begin{x}{\small\bf EXAMPLE} \ (Jordan length for discontinuous curves) \ 
\\[-.5cm]

Let $C : x = x(t)$, $0 \leq t \leq 1$, $x = (x_1, \ldots, x_k)$, 
be any real vector function (a parameteric curve not necessarily continuous).   
Take $A = [0,1]$, $\{I\}$ the class of all intervals $I = [a,b] \subset A$, 
$\sD$ the class of all finite subdivisions 
$D = [I_j, \hsy j = 1, \ldots , N]$, $I_j = [a_{j - 1}, a_j]$, $0 = a_0 < a_1 < \cdots < a_N = 1$.
Let 
\[
s^+ (t_0) 
\ = \ 
\begin{cases}
\limsup\limits_{t \downarrow t_0} \norm{x(t) - x(t_0)} 
\hspace{0.5cm} 0 \leq t_0 < 1
\\[8pt]
\quad 0 \hspace{4cm} t_0 = 1
\end{cases}
.
\]
Analogously, let
\[
s^- (t_0) 
\ = \ 
\begin{cases}
\limsup\limits_{t \uparrow t_0} \norm{x(t) - x(t_0)} 
\hspace{0.5cm} 0 < t_0 \leq 1
\\[8pt]
\quad 0 \hspace{4cm} t_0 = 0
\end{cases}
.
\]
Finally, let
\[
s(t_0) 
\ = \ 
s^+ (t_0) + s^- (t_0) 
\qquad 0 \leq t_0 \leq 1.
\]
As in the continuous situation, suppose that 
\[
\begin{cases}
\ \phi (I) = (\phi_1, \ldots, \phi_k)
\\
\quad \text{where} 
\\
\ \phi_r = x_r (b) - x_r(a) \quad (r = 1, \ldots, k)
\end{cases}
\forall \ I = [a,b] \subset A.  
\]
Then 
\[
V 
\ = \ 
\sup \
\sum \ 
\norm{\phi (I)}
\]
is the Jordan length of $C$, and $V < +\infty$ if and only if $x_r(t)$, $r = 1, \ldots, k$ are BV in $A$.  
If 
$V < +\infty$, 
then $s(t)$, $0 \leq t \leq 1$, is zero everywhere but for countably many $t$, and the sum (or the sum of the series)
\[
\sigma 
\ = \ 
\sum\limits_{0 \leq t \leq 1} \ 
s(t) 
\  \leq \ 
V 
\]
is finite.  
If we take
\[
\delta (D) 
\ = \ 
\max\limits_{i = 1, \ldots, N} \abs{I_i} 
+ 
\sigma 
-
\sum\limits_{i = 0}^N \ 
s(t_i)
\]
where
\[
\begin{cases}
\ D = [I_1, \ldots, , I_N]
\\
\ I_i = [a_{i-1}, a_i]
\quad 
0 = a_0 < a_1 < \cdots < a_N = 1
\end{cases}
,
\]
then $\delta (D)$ is a mesh function and $\phi (I)$, $I \in \{I\}$ is quasi additive w.r.t. $\delta (D)$ and $\sD$.  
We leave the proof to the reader.  
Then it is easy to prove also that 
\begin{align*}
V \ 
&=\ 
V
\left(
\norm{\phi}
\right)
\\[15pt]
&=\ 
\lim\limits_{\delta (D) \ra 0}\ 
\sum\limits_{\substack{I \in D\\D \in \sD}}\ 
\norm{\phi(I)}.
\end{align*}
\\[-1cm]
\end{x}

\begin{x}{\small\bf EXAMPLE} \ (Cauchy integral in an interval in $\R^m$) \ 
\\[-.5cm]

We may suppose $A$ is the unit interval 
\[
A 
\ = \ 
\{(x_1, \ldots, x_m) \in \R^m : 0 \leq x_r \leq 1, r = 1, \ldots, m\}.
\]
Let $\sD$ be the collection of all finite subdivisions $D = I$ of $A$ into intervals 
\[
I 
\ = \ 
[a_r \leq x_r \leq b_r, r = 1, \ldots, m].
\]
Let 
$
f(x), 
\  x \in A 
:
\abs{f(x)} \leq M
$
be any bounded real function.  
For every $I$, let 
$
\begin{cases}
\ \ell = \ell(I) 
\\
\ L = L(I)
\end{cases}
$
denote the \ 
$
\begin{cases}
\ \text{infimum}
\\
\ \text{supremum}
\end{cases}
$
of $f(x)$ in $I$, 
and let
\[
\begin{cases}
\ \Delta (I) = L - \ell
\\
\ F = F(I), \quad F : \ell \leq F \leq L
\\
\ \phi (I) = F(I) \abs{I}, \quad \abs{I} = (b_1 - a_1) \cdots (b_m - a_m)
\end{cases}
.
\]
Let
\[
\delta (D) 
\ = \ 
\max\limits_{\substack{I \in D, \\
D \in \sD}} \ \diam I.
\]
Then $\delta (D)$ is a mesh.  
The Riemann integrability condition reads:
\\[-.25cm]

\qquad (R) \quad Given $\epsilon > 0$ there is a $\sigma = \sigma (\epsilon) > 0$ such that 
\[
\sum \ 
\Delta (I) \abs{I} 
\ < \ 
\epsilon
\quad \forall \ D \in \sD \ \text{with} \ \delta (D) < \sigma.
\]
It is easy to verify that (R) implies quasi additivity of $\phi (I)$ (w.r.t $\sD$ and $\delta (D)$) and vice-versa.  
The limit
\[
\sB 
\ = \ 
\lim \ 
\sum \ 
\phi (I)
\]
is the \un{Cauchy integral} of $f(x)$ in $A$.
\\
\end{x}

\begin{x}{\small\bf EXAMPLE} \ (Lebesgue-Stieltjes integral) \ 
Let $\mu$ be a measure in a $\sigma$-ring $\sR$ of subsets of a space $A$, 
i.e., $(A, \sR, \mu)$ 
is a measure space, and let $f(x)$, $x \in A$, be a $\mu$-measurable and $\mu$-integrable function.  
It is not restrictive to suppose $f(x) \geq 0$.  
For every 
$0 < p < +\infty$, 
and 
$0 \leq p < q \leq +\infty$, 
denote by $B(p)$, $I(p, q)$ respectively the sets 
\[
B(p) 
\ = \ 
\{x \in A : f(x) = p\}
\]
and 
\[
I(p, q)
\ = \ 
\{x \in A : p < f(x) \leq q\}.
\]
These sets are all $\mu$-measurable and 
\[
\begin{cases}
\ 0 \leq p \hsy \mu (B(p)) < +\infty
\\[4pt]
\ \mu(B(+\infty)) = 0
\end{cases}
,  \ \text{and}  \quad p \hsy \mu (B(p)) \equiv 0 \quad \text{for} \quad p = +\infty.
\]
And
\[
\begin{cases}
\ p \hsy  \hsy \mu(I(p, q)) < +\infty
\\[11pt]
\ \ds\sum\limits_{p > 0} \ p \hsy \mu (B(p)) \leq M < +\infty 
\quad (\exists \ M \geq 1) \hspace{1cm}
\end{cases}
.
\hspace{2.25cm}
\]
Thus the set $P$ of all $p > 0$ with $\mu(B(p)) > 0$ is countable and the corresponding series
\[
\sum \ 
p \hsy \mu(B(p))
\]
is convergent.  
Also
\[
\begin{cases}
I (p - \tau, p + \tau) \ra B(p)
\\[8pt]
\mu (I (p - \tau, p + \tau)) \ra \mu (B(p)) 
\end{cases}
\tau \downarrow 0, \ \forall \ p > 0.
\]
Let $\{I\}$ be the collection of all sets 
\[
I(p,q) 
\quad 
0 \leq p < q \leq +\infty,
\]
and let
\[
\psi (I) 
\ = \ 
p \hsy \mu( I(p,q)), 
\quad I \in \{I\}.
\]
Let $\sD$ be the family of all finite decompositions $D = [I_i, \ i = 1, \ldots, n]$ of $A$ into sets
\[
I_i 
\ = \ 
\{x \in A : p_{i - 1} < f(x) \leq p_i\} \in \{I\} 
\]
with
\[
0 = p_0 < p_1 < \cdots < p_{n-1} < p_n = +\infty.
\]
Finally, let 

\[
\delta (D) 
\ = \ 
\begin{cases}
\qquad 1 \hspace{7.5cm} n = 1
\\[15pt]
\ds\ \max\limits_{i = 1, \ldots, n = 1} \ 
(p_i - p_{i-1}) + \frac{1}{p_{n-1}} + 
\sum\limits_{i = 1}^{n-1} \ 
p_i \hsy \mu (B(p))
\quad n > 1
\end{cases}
.
\]
Then $\delta (D)$ is a mesh function and $\psi (I)$ is quasi additive (w.r.t $\sD$ and $\delta (D)$).  
To prove the last statement, given $\varepsilon > 0$, take 
\[
\eta 
\ = \  
\eta(\varepsilon)
\ = \  
\min \hsx \left[\frac{\varepsilon}{4M}, 1\right]
\]
and let
\[
D_0 
\ = \ 
[I_i = I (p_{i - 1}, p_i) : i = 1, \ldots, N], 
\quad 
0 = p_0 < p_1 < \cdots < p_{n-1} < p_n = +\infty
\]
be any $D_0 \in \sD$ with $\delta (D_0) < \eta$. 
Then there are numbers $\tau_i > 0$ such that 
\[
\mu (I(p_i - \tau_i, p_i + \tau_i)
\ < \ 
\mu (B(p_i)) 
+ \frac{\varepsilon}{4 \hsy N \hsy p_i}, \ i = 1, \ldots, N-1.
\]
Take
\[
\tau 
\ = \ 
\min \hsx \tau_i, 
\]
and
\[
\lambda 
\ = \ 
\lambda(\varepsilon, D_0)
\ = \ 
\min \hsx [p_i - p_{i-1}, \hsx i = 1, \ldots, N - 1; 
\hsx 1/p_{N-1};\hsx \tau;\hsx \eta].
\]
If 
\[
D 
\ = \ 
[J_i = (q_{j-1}, q_j), j = 1, \ldots, n], 
\quad 
0 = q_0 < q_1 < \cdots < q_{n-1} < q_n = +\infty,
\]
is any $D \in \sD$ with $\delta (D) < \lambda$, we have 
$J_1 \subset I_1$, $J_n \subset I_N$, and 
\begin{align*}
\Delta \ 
&= \ 
\sum\limits_{i = 1}^N \ 
\left[
\sum{}^{(i)} \ 
\psi (J_i) - \psi (I_i)
\right]^-
\\[15pt]
&= \ 
\sum\limits_{i = 1}^N \ 
p_{i - 1} \hsy 
\left(
\sum{}^{(i)} \ 
\mu (J (q_{j-1}, q_j)) - \mu (I(p_{i-1}, p_i))
\right)^-
\end{align*}
where all $q_{j - 1}$, $q_j$ relative to $\sum{}^{(i)}$ are between $p_{i-1}$ and $p_i$, and $q_0 = p_0 = 0$.  
We have
\begin{align*}
\Delta \ 
&\leq \ 
\sum\limits_{i = 2}^N \ 
p_{i - 1} \hsy 
\left(
\sum{}^{(i)} \ 
\mu (J (q_{j-1}, q_j)) - \mu (I(p_{i-1}, p_i))
\right)^-
\\[15pt]
&\leq \ 
\sum\limits_{i = 2}^N \ 
p_{i - 1} \hsy 
\left[
\mu(B(p_{i-1})) 
+ 
\frac{\varepsilon}{4} N \hsy p_{i-1} 
+ 
\mu(B(p_{i-1})) 
+ 
\frac{\varepsilon}{4 \hsy N \hsy p_i} 
\right]
\\[15pt]
&\hspace{3cm}
+ 
p_{N-1} 
\left[
\mu(B(p_{N-1})) 
+ 
\frac{\varepsilon}{4 \hsy N \hsy p_{N-1}}  
\right]
\\[15pt]
&< \
2 \delta (D_0) 
+ 
\frac{\varepsilon}{4} 
+ 
\frac{\varepsilon}{4} 
\\[15pt]
&< \
\frac{\varepsilon}{2}
+ 
\frac{\varepsilon}{4} 
+ 
\frac{\varepsilon}{4} 
\\[15pt]
&= \ 
\varepsilon.
\end{align*}
Thus $\psi (I)$, $I \in \{I\}$, is quasi subadditive (w.r.t. $\sD$ and $\delta (D)$), 
and hence 
\[
V 
\ = \ 
V(\psi) 
\ = \ 
\lim\limits_{\delta (D) \ra 0} \ 
\sum\limits_{I \in D} \ 
\psi (I).
\]
Since $f(x)$ is $\mu$-integrable, we have 
\[
V 
\ = \ 
V(\psi) 
\ < \ 
+\infty,
\]
and (cf. \S9, \#16), $\psi$ is quasi additive.  
Also it is easy to prove that 
\[
V(\psi) 
\ = \ 
(A) \ 
\int \ 
f(x) 
\ \td \mu
\]
is the Lebesgue-Stieltjes integral of $f(x)$ in $A$ 
(see, e.g. 
Rosenthal\footnote[2]{\textit{Set Functions}, Albuquerque, University of New Mexico Press, (1948).} 
and
Zaanen\footnote[3]{\textit{Theory of Integration}, New York, Interscience, (1958).} 
)  
It is also known that the limit above exists even if taken as 
\[
\lim\limits_{\Delta \ra 0} \ 
\sum\limits_{I \in D} \ 
\psi (I)
\qquad (\Delta \hsx = \hsx \max \ [p_i - p_{i-1}] + 1/p_{n-1})
\]
\\[-1cm]
\end{x}


\begin{x}{\small\bf EXAMPLE} 
(Weierstrass integral over a rectifiable continuous curve $C$)
\\[-.5cm]

Let $C : x = x(t)$, $0 \leq t \leq 1$, $x = (x_1, \ldots, x_k)$, 
be any real continuous vector function with finite Jordan length $L$.  
As in Example \#7, let $A = [0,1]$,  
$I$ the class of all $I = [a,b] \subset A$, $\sD$ the class of all finite subdivisions 
$D = [I_i : i = 1, \ldots N]$, 
$I_i = [a_{i -1}, a_i]$, 
$0 = a_0 < a_1 < \cdots < a_N = 1$, 
and  
\[
\delta 
\ = \  
\delta (D)
\ = \  
\max\hsx \{a_i - a_{i - 1}\}.  
\]
And
\[
\begin{cases}
\ \phi (I) = (\phi_1, \ldots, \phi_k)
\\
\quad \text{where} 
\\
\ \phi_r (I) = x_r (b) - x_r(a) \quad (r = 1, \ldots, k)
\end{cases}
\forall \ I = [a,b] \subset A.  
\]
Then $L = V(\norm{\phi})$ (cf. Example \#7) and the functions
$\phi,  \hsy \phi_r, \hsy \norm{\phi}, \hsy \abs{\phi}, \hsy \phi_r^+, \hsy \phi_r^-$ 
are all quasi additive w.r.t. $D$ and $\delta (D)$.  
\\[-.5cm]

Let $K$ be any compact set containing the graph $[C]$ of $C$, 
i.e., 
\[
[C]
\subset 
K 
\subset
\R^k,
\]
and let 
$f(p, q)$, $p \in K$, $q \in \R^k$, 
be any function continuous on $K \times \R^k$ such that 
\[
f (p, t q) 
\ = \ 
t \hsy f(p, q), 
\quad \forall \  t \geq 0, \hsx p \in K, \hsx q \in \R^k.
\]
For every $I \in \{I\}$, let 
\[
\Phi (I) 
\ = \ 
f [x(\tau), \phi (I)] 
\qquad \tau \in I \quad \text{arbitary}.
\]
Then (\S10, \#7) $\Phi (I)$, $I \in \{I\}$, is quasi additive w.r.t. $\sD$ and $\delta (D)$, 
and the numbers $\eta(\varepsilon)$, $\lambda (\varepsilon, D_0)$ can be taken independently of the 
choice of the $\tau$'s in the intervals $I$.  
Thus, the following limit exists and is finite:
\begin{align*}
\sF \ 
&= \ 
\int\limits_C \ 
f(p, q) 
\\[15pt]
&= \ 
\lim\limits_{\delta (D) \ra 0}\ 
\sum\limits_{I \in D} \ 
\Phi (I)
\\[15pt]
&= \ 
\lim\limits_{\delta (D) \ra 0}\ 
\sum\limits_{i = 1}^n \ 
f[x(\tau_i), x(a_i) - x(a_{i - 1})],
\end{align*}
independently of the choices of the points $\tau_i \in [a_{i-1}, a_i]$.  
This limit is known as the \un{Weierstrass integral} of $f$ on the rectifiable curve $C$ 
and was studied by 
L. Tonelli, 
G. Bouligand, 
N. Aronszajn, 
K. Menger, 
C. Y. Pauc.   
See
Cesari\footnote[2]{\textit{Surface Area}, Annals of Mathematics Studies, Number \textbf{35}, Princeton University Press, (1956).} 
for references.
\\[-.5cm]

It is known that $\sF$ is invariant w.r.t. Fr\'echet equivalence, i.e., it is independent of the representation $C : x = x(t)$ of $C$.  
In particular, if $s = s(t)$ denotes the Jordan length of the curve $C_t$ defined by $x(t)$ on $[0,t]$, $0 \leq s \leq L$, 
and $X (s) = x [s(t)]$, then $C : x = x(s)$, $0 \leq s \leq L$, is a representation of $C$, $X(s)$ is Lipschitz of constant 1, 
and $\sF$ is given by the Lebesgue integral
\begin{align*}
\sF \ 
&= \ 
\int\limits_C \ 
f(p, q) 
\\[15pt]
&= \ 
\int\limits_0^L \ 
f[X(s), X^\prime (s)] 
\ \td s.
\end{align*}
\\[-1cm]
\end{x}


\centerline{\textbf{\large REFERENCES}}
\setcounter{page}{1}
\setcounter{theoremn}{0}
\renewcommand{\thepage}{References-\arabic{page}}
\vspace{0.75cm}

\[
\text{BOOKS}
\]

\begin{rf}
Cesari, Lamberto, \textit{Surface Area}, Annals of Mathematics Studies, Number 35, Princeton University Press, 1956.
\end{rf}

\begin{rf}
Evans, L. C. and Gariepy, R. F., \textit{Measure theory and fine properties of functions}, CRC Press, Boca Raton, 1992.
\end{rf}

\begin{rf}
Goffman, Casper, Nishiura, Toga, and Waterman, Daniel, \textit{Homeomorphisms in Analysis}, American Mathematical Society, 1997.
\end{rf}

\begin{rf}
Goffman, C., \textit{Real Functions}, Rinehart, New York, 1953.
\end{rf}

\begin{rf}
Leoni, Giovanni, \textit{A First Course in Sobolev Spaces}, Second Edition, American Mathematical Society, Providence, Rhode Island, 2017.  (pp. 133-155).
\end{rf}


\setcounter{theoremn}{0}

\[
\text{ARTICLES}
\]

\begin{rf}
Burkill, J.C., Functions of Intervals, \textit{Proc. London Math.} Soc. (2) vol. \textbf{22} (1924), 275-310.
\end{rf}

\begin{rf}
Cesari, L., Quasi-additive set functions and the concept of integral over a variety, \textit{Trans. Amer. Math. Soc.} vol. \textbf{102} (1962), 94-113.
\end{rf}

\begin{rf}
Cesari, L., Variation, Multiplicity, and Semicontinuity, \textit{American Mathematical Monthly} vol. \textbf{65}, issue 5, (1958), 317-332.
\end{rf}

\begin{rf}
Cesari, L., Rectifiable Curves and the Weierstrass Integral, \textit{American Mathematical Monthly} vol. \textbf{65}, issue 7, (1958), 485-500.
\end{rf}

\begin{rf}
Fr\'echet, M., Sur l’integral d’une fonctionnelle etendue a un ensemble abstrait, \textit{Bull. Soc. Math.} France. vol. \textbf{43} (1915), 249-267.
\end{rf}

\begin{rf}
Gavurin, M., Uber die Stieltjessche Integration abstrakter Funktionen, \textit{Fund. Math.} vol. \textbf{27} (1936), 255-268.
\end{rf}

\begin{rf}
Goffman, C., Non-parametric surfaces given by linearly continuous functions, \textit{Acta Math.} \textbf{103} (1960), 269-291.
\end{rf}

\begin{rf}
Goffman, C., A characterization of linearly continuous functions whose partial derivatives are measures, \textit{Acta Math.} \textbf{117} (1967), 165-190.
\end{rf}

\begin{rf}
Henstock, R., On interval functions and their integrals, \textit{J. Lond. Math. Soc.} vol. \textbf{21} (1946), 204-209.
\end{rf}

\begin{rf}
Kempisty, S., Fonctions d’intervalle non additives, \textit{Actualites Sci. Ind.} No. \textbf{824}, Paris Hermann, 1939.
\end{rf}

\begin{rf}
Kober, H., On the existence of the Burkill integral, \textit{Canad. J. Math.} vol. \textbf{10} (1957), 115-121.
\end{rf}

\begin{rf}
Kolmogoroff, A., Untersuchungen uber den Intergralbergriff, \textit{Math. Ann.} vol. \textbf{103} (1930), 654-696.
\end{rf}

\begin{rf}
Krickeberg, K., Distributionen, Funktionen beschrankter variation und Lebesguescher Inhalt nichtparametrischer Flachen, \textit{Ann. Mat. Pura Appl.} \textbf{4} (1957), 105-133.
\end{rf}

\begin{rf}
McShane, E. J., Generalized Curves, \textit{Duke Math. Journal} vol. \textbf{6} (1940), 513-536.
\end{rf}

\begin{rf}
Pollard, S., The stieltjes integral and its generalizations, \textit{Quart. J. Math. vol.} \textbf{49} (1923), 73-138.
\end{rf}

\begin{rf}
Serrin, J., A new definition of the integral for non-parametric problems in the calculus of variations, \textit{Acta Math.} vol. \textbf{102} (1959), 23-32.
\end{rf}

\begin{rf}
Smith, H. L., Stieltjes Integral, \textit{Trans. Amer. Math. Soc.} vol. \textbf{27} (1925), 490-507.
\end{rf}

\begin{rf}
Warner, G., The Burkill-Cesari Integral, \textit{Duke Mathematical Journal} vol. \textbf{35}, \#1, (1968), 61-78.
\end{rf}

\newpage

\setcounter{theoremn}{0}

\[
\text{Ph. D. THESES}
\]

\begin{rf}
Breckenridge, John Clark, Geocze K-Area and Measure Theoretical Methods in Surface Area Theory, 1969.
\end{rf}

\begin{rf}
Gariepy, Ronald Francis, Current Valued Measures and Geocze Area, 1969.
\end{rf}

\begin{rf}
Nishiura, T., Analytic Theory of Continuous Transformations, 1959.
\end{rf}

\begin{rf}
Shapley, L. S., Additive and Non-Additive Set Functions, 1953.
\end{rf}

\begin{rf}
Stoddart, A. W. J., Integrals of the Calculus of Variations, 1964.
\end{rf}

\begin{rf}
Turner, L. H., The Direct Method in the Calculus of Variations, 1957.
\end{rf}

\begin{rf}
Warner, Garth W., Quasi Additive Set Functions and Non-linear Integration over a Variety, 1966.
\end{rf}%
\printindex
\end{document}